\documentclass[12pt,reqno,a4paper]{amsart}
\usepackage{amsmath}
\usepackage{amssymb}
\usepackage{stmaryrd}
\usepackage{amsthm}

\usepackage[applemac]{inputenc}

\def\EM#1{{\em #1\/}}
\newif\ifenglish
\englishtrue

\newtheorem{thm}{\ifenglish Theorem\else Satz\fi}
\newtheorem{pro}{Proposition}
\newtheorem{lem}{Lemma}

\newtheorem{cor}{\ifenglish Corollary\else Korollar\fi}

\newtheorem{obs}{\ifenglish Observation\else Beobachtung\fi}

\def\figref#1{Figure~\ref{#1}}

\def\defeq{:=}

\def\setof#1{\left\{#1\right\}}
\def\of#1{\!\left(#1\right)}

\def\pas#1{\left(#1\right)}
\def\brk#1{\left[#1\right]}
\def\sgn{{\operatorname{sgn}}}

\def\floor#1{{\left\lfloor #1\right\rfloor}}
\def\ceil#1{{\left\lceil #1\right\rceil}}
\def\bit{\begin{itemize}}
\def\eit{\end{itemize}}
\def\beq{\begin{equation}}
\def\eeq{\end{equation}}

\def\N{{\mathbb N}}

\def\Z{{\mathbb Z}}

\def\1{{\mathbf 1}}

{\begin{list}{#1}{\parsep0em \itemsep0.3em \labelwidth1em
\labelsep0.5em \leftmargin1.5em }} {\end{list}}

\def\lgt{Lindstr\"om--Gessel--Viennot}

\def\detof#1{\left\vert{#1}\right\vert}

\def\path{{%\mathrm 
P}}
\def\figref#1{Figure~\ref{#1}}

\def\symfs{\Lambda}	
\def\csymf{h}		

\def\weight{\omega}

\def\complement#1{\overline{#1}}

\def\cardof#1{\left\vert{#1}\right\vert}                

\def\symm{\mathfrak{S}}

\def\and{\wedge}

\def\complement#1{\overline{#1}}

\def\myrange#1{\setof{1,\dots,k}}

\def\schurf{s}
\def\length#1{\ell\of{#1}}
\def\tableau{T}

\usepackage{pstricks}
\usepackage{multido}
\usepackage{graphicx}
\usepackage{color}
\newcommand{\showgrid}{}
\newcommand{\gridon}{\renewcommand{\showgrid}{\psset{subgriddiv=1,griddots=10,gridlabels=6pt}\psgrid}}
\gridon

\advance\voffset    by -1  cm
\advance\hoffset    by -1.5cm
\advance\textwidth  by  3  cm
\advance\textheight by  1  cm

\sffamily

\def\n{n}		

\def\cop{{\sl\bf cop}}
\def\coc{{\sl\bf rop}}
\def\nop{{\sl\bf nop}}
\def\pointconf{\EM{{\bf c}onfiguration {\bf o}f (lower/upper) {\bf p}oints} (short: \cop){\gdef\pointconf{\cop}}}
\def\pointarr{\EM{{\bf r}adial orientation {\bf o}f (coloured lower/upper) {\bf p}oints} (short: \coc){\gdef\pointarr{\coc}}}
\def\pointmatch{\EM{{\bf n}oncrossing {\bf o}riented {\bf p}erfect matching} (short: \nop){\gdef\pointmatch{\nop}}}
\def\nilp{non\-in\-ter\-sec\-ting lat\-ti\-ce path}

\def\mdetof#1{\left\vert#1\right\vert}

\def\jtmat#1{{\mathbf h}_{#1}} 
\def\deleterowcols#1#2#3{\brk{#1}_{\underline{#2}\vert\underline{#3}}}
\def\deleterowcols#1#2#3{\brk{#1}_{\text{\sout{$#2$}}\vert\text{\sout{$#3$}}}}
\def\deleterow#1#2{\deleterowcols{#1}{#2}{\hphantom{0}}}
\def\deletecol#1#2{\deleterowcols{#1}{\hphantom{0}}{#2}}
\def\range#1{\brk{#1}}
\def\shiftrange#1#2{\pas{\range#1+#2}}
\def\sumset#1{{\scriptstyle\pas{\sum_{x\in #1}x}}}
\def\sumset#1#2{\Sigma_{#1\subseteq#2}}
\def\shuffle#1#2#3{\pas{{\scriptscriptstyle #2\!\leftarrow}\!#1\!{\scriptscriptstyle\leftarrow \!#3}}}
\def\shuffle#1#2#3{\pas{\left.#1\right\vert^{{\scriptscriptstyle#2\rightarrow}}_{{\scriptscriptstyle\leftarrow#3}}}}
\def\minor#1#2#3{\pas{#1}_{#2,\;#3}}
\def\deleterowcols#1#2#3{\minor{#1}{\overline{#2}}{\overline{#3}}}
\def\deleterow#1#2{\deleterowcols{#1}{\setof{#2}}{\emptyset}}
\def\deletecol#1#2{\deleterowcols{#1}{\emptyset}{\setof{#2}}}

\def\detprodcminors#1#2#3{\detof{\minor{#1}{#2}{#3}}\cdot\detof{\deleterowcols{#1}{#2}{#3}}}
\def\limes#1{#1_\infty}
\def\allconstpart#1{\pas{{#1}}}
\def\constpart#1#2{\pas{{#1}^{\pas{#2}}}}
\def\shape#1#2{#1/#2}
\def\partweightof#1{\cardof{#1}}
\def\deletepart#1#2{{#1}^{\overline{#2}}}
\def\addtoparts#1#2{\pas{#1+\constpart{#2}{m}}}

\def\harc#1{{\scriptstyle #1}^\to}
\def\harc#1{\overline{#1}}
\def\harc#1{a^h_{#1}}
\def\varc#1{{\scriptstyle #1}^\uparrow}
\def\varc#1{{#1}\vert}
\def\varc#1{a^v_{#1}}
\def\nofentries#1#2{\#\of{#1,#2}}
\def\gfof#1{{\mathbf{gf}\of{#1}}}
\def\paths#1#2{{\mathfrak P}\of{#1\to#2}}
\def\pathsmod#1#2#3{{\mathfrak P}^{#3}\of{#1\to#2}}
\def\gfpaths#1#2{\gfof{\paths{#1}{#2}}}
\def\idperm{\operatorname{id}}
\def\alljt{{\mathfrak D}}
\def\ijt{{\mathfrak I}}
\def\nijt{{\mathfrak N}}
\def\involution#1{{\mathbf{#1}}}
\def\invol{\involution{i}}
\def\asvector#1{{\mathbf #1}}
\def\subpart{\trianglelefteq}
\def\identity#1{{\mathcal #1}}
\def\cnr#1{{\mathbf #1}}

\newcounter{dummycount}
\def\TRIup#1{\pspolygon[linewidth=0.1,linecolor=black,fillstyle=solid,fillcolor=#1](0.0,1.25)(-0.33,0.75)(0.33,0.75)}
\def\TRIdown#1{\pspolygon[linewidth=0.1,linecolor=black,fillstyle=solid,fillcolor=#1](0.0,0.75)(-0.33,1.25)(0.33,1.25)}

\newif\iflongversion
\longversionfalse

\begin{document}

\bibliographystyle{plain}

\title[Viewing determinants as nonintersecting lattice paths]{Viewing determinants as nonintersecting lattice paths yields
	classical determinantal identities bijectively}

\begin{abstract}
In this paper, we show how general determinants may be viewed as generating functions of
nonintersecting lattice paths, using the \lgt--interpretation of semistandard
Young tableaux and the Jacobi--Trudi identity together with elementary observations.
After some preparations, this point of view provides very simple ``graphical proofs''
for classical determinantal identities like the Cauchy--Binet formula, Dodgson's
condensation formula, the Pl\"ucker relations and Laplace's expansion.
Also, a determinantal identity generalizing Dodgson's condensation formula is presented,
which might be new.
\end{abstract}

\author{Markus Fulmek}
\address{Fakult\"at f\"ur Mathematik, 
Nordbergstra\ss e 15, A-1090 Wien, Austria}
\email{{\tt Markus.Fulmek@Univie.Ac.At}\newline\leavevmode\indent
{\it WWW}: {\tt http://www.mat.univie.ac.at/\~{}mfulmek}
}

\date{\today}
\thanks{
Research supported by the National Research Network ``Analytic
Combinatorics and Probabilistic Number Theory'', funded by the
Austrian Science Foundation. 
}

\maketitle

\section{Introduction}\label{sec:intro}
In \cite{fulmek:2001}, a combinatorial proof was given for two Schur function
identities, which were presented in \cite{kirillov:1984} and in \cite{kleber:2001}.
This combinatorial proof was shown to apply to a \EM{class} of Schur function identities
\cite[Lemma~16]{fulmek:2001}, and was used to prove bijectively Dodgson's condensation
formula and the Pl\"ucker relations for examples, but was not paid further attention.
Recently, members of this class of Schur function identities received some interest
\cite{gurevichPyatovSaponov:2009}. The close connection of the result \cite[(3.3)]{gurevichPyatovSaponov:2009} was already explained \EM{adhoc} in \cite{fulmek:2009},
but  we take this opportunity to make obvious the much wider range of applicability of
this idea, which amounts to ``viewing determinants as (generating functions of) \nilp s'',
by giving concrete examples. The combinatorial constructions
are best conceived by pictures, so we give a lot of illustrations.

This paper is organized as follows: 

In Section~\ref{sec:background}, we present basic background information regarding
symmetric functions, partitions, Young tableaux and (skew) Schur functions.

In Section~\ref{sec:main}, we present the \lgt--interpretation of semistandard
Yount tableaux as \nilp s, and illustrate this view by giving a ``graphical proof''
of the Cauchy--Binet formula.

In Section~\ref{sec:overlays}, we present the central bijective construction
(recolouring of bicoloured trails in the overlays of families
of \nilp s corresponding to some product of skew
Schur functions) and indicate how this construction applies to a class of Schur function
identities.

In Section~\ref{sec:apps}, we present several examples: We prove 
a generalization of Dodgson's condensation rule which might be new (Theorem~\ref{thm:maybenew})
and give
``graphical proofs'' for the Pl\"ucker relations (and its generalization
\cite[(3.3)]{gurevichPyatovSaponov:2009}) and for (a generalization of)
Laplace's expansion.

\section{Basic definitions}\label{sec:background}
The notation $\detof{x}$ has three different meanings in our presentation, depending on the
type of object $x$ (the respective meaning should always be clear from the context):
\bit
\item if $x$ is a \EM{set}, then $\cardof{x}$ denotes the \EM{cardinality} of $x$,
\item if $x$ is a \EM{matrix}, then $\detof{x}$ denotes the \EM{determinant} of $x$,
\item if $x$ is a \EM{partition} or \EM{shape}, then $\partweightof{x}$ denotes the
	\EM{sum of parts} of $x$ (to be explained below).
\eit

In the following, we shall briefly recall basic concepts and facts needed for our presentation. (More information can be found, e.g., in \cite{macdonald:1995}.)

\subsection{The ring of symmetric functions}
Consider the ring $\Z\brk{x_1,x_2,\dots,x_\n}$ of polynomials in $\n$ independent variables
$\asvector{x}\defeq\pas{x_1,x_2,\dots,x_\n}$ with integer coefficients. 
The \EM{degree} of a monomial $x_1^{k_1}\cdot x_2^{k_2}\cdots x_\n^{k_\n}$ is
the sum $k_1+k_2+\dots+k_\n$, and a polynomial $p$ is called \EM{homogeneous} of degree $k$
if all monomials of $p$ have the same degree $k$.

The symmetric group $\symm_\n$ acts on this ring by permuting the variables, and a polynomial
is \EM{symmetric} if it is invariant under this action. The set of all symmetric polynomials
forms a subring $\symfs_n\subseteq\Z\brk{x_1,x_2,\dots,x_\n}$ which is \EM{graded}, i.e.,
\begin{equation*}
\symfs_\n=\bigoplus_{k\geq 0}\symfs_\n^k,
\end{equation*}
where $\symfs_\n^k$ consists of the homogeneous symmetric polynomials of degree $k$, together
with the zero polynomial.

For each $r\in\Z$, the \EM{complete symmetric function} $\csymf_r\of{\asvector{x}}$ is the sum of
all monomials of degree $r$. 
In particular,
$\csymf_0\of{\asvector{x}}=1$ and, by convention, $\csymf_r\of{\asvector{x}}=0$ for $r<0$. 

For example, if $\n=3$, then $\asvector{x}=\pas{x_1,x_2,x_3}$ and
$$
\csymf_2\of{x_1,x_2,x_3}= x_1^2 + x_2^2 + x_3^2 + x_1 x_2 + x_1 x_2 + x_2 x_3.
$$

The following fact is well--known (see, e.g., \cite[(2.8)]{macdonald:1995}):
\begin{pro}
\label{pro:hsym-alg-ind}
The set of all homogeneous symmetric functions is algebraically independent, i.e., $p\equiv0$
is the \EM{only} polynomial such that $p\of{\csymf_0\of{\asvector{x}},\csymf_1\of{\asvector{x}},\dots} \equiv 0$.

Moreover, for $0<k<\n$, the set
$$ \setof{\csymf_{i}\of{x_1,\dots,x_k}:\;i=0,1,2,\dots}
\cup
\setof{\csymf_{i}\of{x_{k+1},\dots,x_{\n}}:\;i=1,2,3,\dots}
$$
is also algebraically independent.
\hfill\qedsymbol
\end{pro}

\subsection{Schur functions}
Recall the following standard definitions:
An infinite weakly decreasing series of nonnegative integers $\lambda=\pas{\lambda_i}_{i=1}^\infty$,
where only finitely many elements are positive, is called a \EM{partition}. The largest
index $i$ for which $\lambda_i>0$ is called the \EM{length} of the partition $\lambda$ and
is denoted by $\length{\lambda}$. The sum of the non--zero parts
$\lambda_1+\lambda_2+\cdots+\lambda_{\length{\lambda}}$ of $\lambda$ is denoted by
$\partweightof\lambda$. In most cases we shall omit the trailing
zeroes, i.e.,
for $\length{\lambda}=r$ we simply write $\lambda=\pas{\lambda_1,\lambda_2,\dots,\lambda_r}$,
where $\lambda_1\geq\lambda_2\geq\dots\geq\lambda_r> 0$.

For example, $\lambda=\pas{7, 4, 4, 3, 1, 1, 1}$ is a partition of length $\length{\lambda}=7$
with $\cardof\lambda=21$.

\begin{figure}
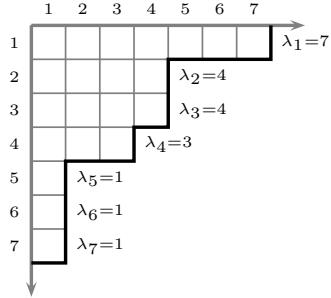

\caption{Illustration: Ferrers diagram $F_\lambda$ of the partition
$\lambda=\pas{7, 4, 4, 3, 1, 1, 1}$.}
\label{fig:FerrersDiagram}
\input graphics/FerrersDiagram-slide
\end{figure}

The {\em Ferrers diagram\/} $F_\lambda$
of $\lambda$
is an array of cells with $\length{\lambda}$ left-justified rows and $\lambda_i$
cells in row $i$. For an illustration, see \figref{fig:FerrersDiagram}.

For our
purposes, it is convenient to generalize this definition: By a \EM{semipartition} we understand
an infinite weakly decreasing series of integers $\pas{\lambda_i}_{i=1}^\infty$,
where
$$\limes\lambda\defeq\lim_{n\to\infty}\lambda_n>-\infty.$$
The length of semipartition $\lambda$ is the largest integer $m$ with
$\lambda_m>\lambda_\infty$, which we denote again by $\length\lambda$:
Note that every partition $\mu$ is a semipartition with $\limes\mu=0$.

Clearly, the set of semipartitions is closed under component--wise addition
$$ \lambda+\mu=\pas{\lambda_i}_{i=1}^\infty + \pas{\mu_i}_{i=1}^\infty\defeq
	\pas{\lambda_i+\mu_i}_{i=1}^\infty.$$
For $m,z\in\Z$, $m>0$, denote by $\allconstpart{z}$ and $\constpart{z}{m}$ the
semipartitions
$$\allconstpart{z}\defeq\pas{z}_{i=1}^\infty\text{ and }\constpart{z}{m}\defeq\pas{\underbrace{z,z,\dots,z}_{m\text{ times}},0,0,\dots},$$
respectively. (Note that $\limes{\constpart{z}{m}}=0$.)
If two semipartitions
$\lambda$, $\mu$ satisfy
\bit
\item $\mu_i\leq\lambda_i$ for all $i=1,2,\dots$,
\item $\limes{\mu} = \limes{\lambda}$,
\eit
then we denote this by $\mu\subpart\lambda$ and introduce the symbol $\lambda/\mu$,
which we call a \EM{shape}. The length of the shape $\lambda/\mu$ is defined by
$\length{\lambda/\mu}\defeq\length\lambda$,
and the (terminating!) sum 
$\sum_{i=1}^\infty\pas{\lambda_i-\mu_i}$ is denoted by
$\partweightof{\lambda/\mu}$. Note that we may view \EM{partition} $\lambda$ as the
\EM{shape} $\lambda/\allconstpart{0}$

\begin{figure}
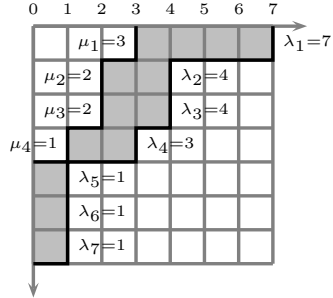

\caption{Illustration: Ferrers diagram $F_{\lambda/\mu}$ of the shape
$\lambda/\mu=\pas{8, 5, 5, 4, 2, 2, 2}/\pas{3, 2, 2, 1}$. The
\EM{same} Ferrers diagram would arise for the shapes $\pas{\lambda+\allconstpart{z}}/\pas{\mu+\allconstpart{z}}$ and $\pas{\lambda+\constpart{z}{7}}/\pas{\mu+\constpart{z}{7}}$, for arbitrary $z\in\Z$.}
\label{fig:skewFerrersDiagram}
\input graphics/skewFerrersDiagram
\end{figure}

The {\em Ferrers diagram\/} $F_{\lambda/\mu}$ of shape $\lambda/\mu$
is an array of cells with $\length{\lambda/\mu}$ left--justified rows and $\pas{\lambda_i-\mu_i}$
cells in row $i$, where the first $\mu_i$ cells in row $i$ are missing, see \figref{fig:skewFerrersDiagram} for an illustration.

Note that for arbitrary $z\in\Z$, the Ferrers diagram $F_{\lambda/\mu}$ also
can be written as
\begin{equation}
\label{eq:FerrersShiftInvariance}
F_{\lambda/\mu} = F_{\lambda+\allconstpart{z}/\mu+\allconstpart{z}}
 = F_{\lambda+\constpart{z}{\length{\lambda/\mu}}/\mu+\constpart{z}{\length{\lambda/\mu}}}.
\end{equation}
In particular,
$F_\lambda$ of partition $\lambda$ may be viewed as the
Ferrers diagram $F_{\lambda+\constpart{z}{\infty}/\constpart{z}{\infty}}$ for \EM{arbitrary} $z\in\Z$. 


\EM{Schur functions}, which are irreducible general linear characters, can be defined
as quotient of alternants \cite{schur:1901} as follows. Let $\lambda$ be a partition and let $\setof{x_1, \dots, x_\n}$ be a set of independent variables.
Then the Schur function $\schurf_\lambda\of{x_1,\dots, x_\n}$ indexed by $\lambda$ is defined
as the quotient of determinants  (see \cite[(3.1)]{macdonald:1995})
\begin{equation}
\label{eq:schur-determinant}
\schurf_\lambda\of{x_1,\dots, x_\n}\defeq
\frac{
	\detof{x_i^{\lambda_j+\n-j}}_{i,j=1}^\n
}{
	\detof{x_i^{\n-j}}_{i,j=1}^\n
}.
\end{equation}
It is easy to see that $\schurf_\lambda$ is a \EM{symmetric function}, which is
homogeneous of degree $\cardof\lambda$. The Jacobi--Trudi identity
(first obtained by Jacobi \cite{jacobi:1841} and simplified by Trudi \cite{trudi:1864},
see \cite[(3.4)]{macdonald:1995}) states that the
Schur function $\schurf_\lambda$ equals the following determinant of complete homogeneous functions:
\begin{equation}
\label{eq:JacobiTrudi}
\schurf_\lambda=
\detof{\csymf_{\lambda_j-j+i}}_{i,j=1}^{\length{\lambda}}.
\end{equation}
Here, we introduced in passing the shortened notations $\schurf_\lambda$ and $\csymf_r$ 
for $\schurf_\lambda\of{x_1,\dots,x_\n}$ and $\csymf_r\of{x_1,\dots,x_\n}$,
respectively.

An {\em $\n$--semistandard Young tableau\/} of shape $\lambda$ is a filling of the cells
of the Ferrers diagram $F_\lambda$ with integers from the set $\setof{1,2,\dots,\n}$, such
that the numbers filled into the cells weakly increase in rows and strictly increase in columns.

\begin{figure}
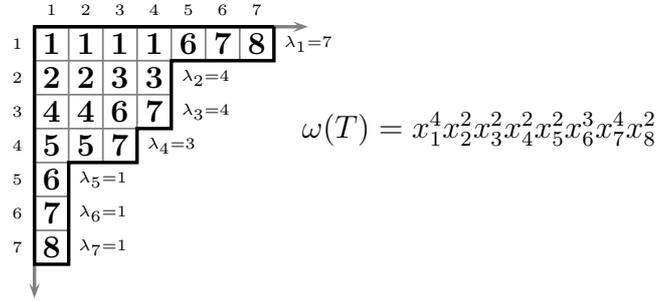

\caption{Illustration: An $8$--semistandard Young tableau $T$ of shape
$\lambda=\pas{7, 4, 4, 3, 1, 1, 1}$ and its weight $\weight\of T$.}
\label{fig:SSYT}
\input graphics/SSYT
\end{figure}

Let ${\tableau}$ be a semistandard Young tableau and define $\nofentries{\tableau}{k}$ to be
the number of entries $k$ in ${\tableau}$. The \EM{weight} $\weight\of{\tableau}$ of ${\tableau}$
is defined as follows:
\begin{equation}
\label{eq:dfn-weight-tableau}
\weight\of{\tableau} = \prod_{k = 1}^{\n} x_k^{\nofentries{\tableau}{k}}.
\end{equation}
See \figref{fig:SSYT} for an illustration.

Then the \EM{Schur function} $\schurf_\lambda$ can equivalently be written as the following
\EM{generating function} (formal sum of weights)
\begin{equation*}
\schurf_\lambda
= \sum_{{\tableau}}\weight\of{\tableau},
\end{equation*}
where the sum is over all $\n$--semistandard Young tableaux ${\tableau}$ of
shape $\lambda$ (see \cite[Definition~4.4.1]{sagan:2000}). 


An {\em $\n$--semistandard skew Young tableau\/} of shape $\lambda/\mu$ is a filling of the cells of $F_{\lambda/\mu}$ with integers from the set
$\{1,2,\dots,\n\}$,
such that the numbers filled into the cells weakly increase in rows and
strictly increase in columns.
See the left picture in \figref{fig:skew-Young-tableau} (or \figref{fig:SSYT} again)
for an illustration.

Then we can define the \EM{skew Schur function} $\schurf_{\lambda/\mu}$
as the following generating function:
\begin{equation}
\label{eq:skewSchur}
\schurf_{\lambda/\mu}
\defeq \sum_{{\tableau}}\weight\of{\tableau},
\end{equation}
where the sum is over all $\n$--semistandard skew Young tableaux ${\tableau}$ of
shape $\lambda/\mu$, where the weight $\weight\of{\tableau}$ of ${\tableau}$
is defined as in \eqref{eq:dfn-weight-tableau}. $\schurf_{\lambda/\mu}$ also is
a symmetric function (see \cite[proof of Proposition 4.4.2]{sagan:2000}), which is
homogeneous of degree $\partweightof{\shape\lambda\mu}$.
%
%
By \eqref{eq:FerrersShiftInvariance}, we clearly have for arbitrary $z\in\Z$:
\begin{equation}
\label{eq:SchurShiftInvariance}
\schurf_{\lambda/\mu} = \schurf_{\lambda+\allconstpart{z}/\mu+\allconstpart{z}}
 = \schurf_{\lambda+\constpart{z}{\length{\lambda/\mu}}/\mu+\constpart{z}{\length{\lambda/\mu}}}.
\end{equation}
%
In particular,
$\schurf_{\lambda+\constpart{z}{\infty}/\constpart{z}{\infty}}
= \schurf_{\lambda+\constpart{z}{\length{\lambda/\mu}}/\constpart{z}{\length{\lambda/\mu}}}$
is identical to the ``ordinary'' Schur function $\schurf_\lambda$.

Skew Schur functions, too, have an expansion in terms of complete homogeneous functions (see \cite[(5.4)]{macdonald:1995}) which generalizes \eqref{eq:JacobiTrudi}:
\begin{equation}
\label{eq:JacobiTrudiSkew}
\schurf_{\lambda/\mu}=
\detof{\csymf_{\lambda_j-\mu_i-j+i}}_{i,j=1}^{\length{\lambda}}.
\end{equation}

\subsection{Connection between determinantal relations and Schur function identities}
We want to illustrate the connection between determinants and skew Schur functions,
using Dodgson's condensation formula as an example.

Consider some $m\times n$--matrix $a=\pas{a_{i,j}}_{\pas{i,j}=\pas{1,1}}^{\pas{m,n}}$ and define $\range{k}\defeq \setof{1,2,\dots,k}$ for $k\in\N$. Then we may write
$a=\minor{a}{\range{m}}{\range{n}}\defeq\pas{a_{i,j}}_{\pas{i,j}\in\range{m}\times\range{n}}$.
More generally, for subsets $R\subseteq\brk{m}$ and $C\subseteq\brk{n}$, we denote
by $\minor{a}{R}{C}$ the \EM{minor} of $a$ which consists of the rows $R$ and the columns $C$;
\EM{in the same order as in $a$.} (All sets we consider in this paper are
\EM{ordered}, and all subsets ``inherit the order''.)

If we want to describe the same minor by deleting rows and columns in $a$, then we
write
$$
\minor{a}{R}{C} =
\deleterowcols{a}{X}{Y},
$$
where $X=\complement{R}\defeq\range{m}\setminus R$ and
$Y=\complement{C}\defeq\range{n}\setminus C$.


\begin{pro}[Dodgson's Condensation] Let $a=\pas{a_{i,j}}_{i,j=1}^m$ be an arbitrary
$m\times m$--matrix, $m\geq2$. Then there holds the following identity:
\begin{equation}
\label{eq:Dodgson}
\mdetof{a}\cdot\mdetof{\deleterowcols{a}{\setof{1,m}}{\setof{1,m}}} =
	\mdetof{\deleterowcols{a}{\setof1}{\setof1}}\cdot\mdetof{\deleterowcols{a}{\setof m}{\setof m}} - 
\mdetof{\deleterowcols{a}{\setof1}{\setof m}}\cdot\mdetof{\deleterowcols{a}{\setof m}{\setof1}}.
\end{equation}
(Note that for $m=2$, this amounts to the formula for $2\times2$--determinants.)\hfill\qedsymbol
\end{pro}

%

We shall ``translate'' this determinantal identity to a Schur function identity. To this end,
we introduce the following operation ``delete part
$\lambda_k$
in semipartition $\lambda$'':
\begin{equation}
\label{eq:deletepart}
\deletepart{\lambda}{k}\defeq\pas{\lambda_1,\dots,\lambda_{k-1},\lambda_{k+1}-1,\lambda_{k+2}-1,\dots}.
\end{equation}
Note that
$\limes{\deletepart{\lambda}{k}}=\limes\lambda-1$
and
$\length{\deletepart{\lambda}{k}}=\max\of{\length\lambda,k}-1$.

For some subset $\setof{k_1<k_2<\dots<k_l}\subseteq
\setof{1,2,\dots,m}$, we define inductively:
$$
\deletepart{\lambda}{k_1,\dots,k_l}\defeq
\deletepart{\pas{\deletepart{\lambda}{k_2,\dots,k_{l}}}}{k_1}.
$$

Now assume some skew shape $\lambda/\mu$, $m=\length{\lambda/\mu}$. According
to \eqref{eq:JacobiTrudiSkew}, the skew Schur function $\schurf_{\lambda/\mu}$ equals the
determinant of the matrix
$$\jtmat{\lambda/\mu}\defeq\pas{\csymf_{\lambda_j-\mu_i-j+i}}_{\pas{i,j}={1,1}}^{\pas{m,m}}.$$
Observe that 
the $\pas{i,j}$--entries in matrix $\deletecol{\jtmat{\lambda/\mu}}{k}$ are
\bit
\item $\csymf_{\lambda_{j}-\mu_i-j+i}$ for $j<k$,
\item $\csymf_{(\lambda_{j+1}-1)-\mu_i-j+i}$ for $k\leq j\leq m-1$
\eit
(i.e., deleting column $k$ in $\jtmat{\lambda/\mu}$ corresponds to deleting $\lambda_k$ in $\lambda$),
while the $\pas{i,j}$--entries in matrix $\deleterow{a_{\lambda/\mu}}{l}$ are
\bit
\item $\csymf_{\lambda_{j}-\mu_i-j+i}$ for $i<l$,
\item $\csymf_{\lambda_{j}-\pas{\mu_{i+1}-1}-j+i}$ for $l\leq i\leq m-1$
\eit
(i.e., deleting row $l$ in $\jtmat{\lambda/\mu}$ corresponds to deleting $\mu_l$ in $\mu$).
These observations generalize to the following relation:
\begin{equation}
\label{eq:minorspartitions}
\deleterowcols{\jtmat{\lambda/\mu}}{\setof{i_1,\dots,i_k}}{\setof{j_1,\dots,j_l}} = \jtmat{\deletepart{\lambda}{i_1,\dots,i_k}/\deletepart{\mu}{j_1,\dots,j_l}}.
\end{equation}
%
So we can deduce the following Schur function identity:
\begin{cor}
\label{cor:Dodgson}
Let $\lambda=\pas{\lambda_1,\lambda_2,\dots,\lambda_m}$ be a partition of length $\length{\lambda}=m\leq2$.
Then we have the following Schur function identiy:
\begin{multline}
\label{eq:DodgsonSchur}
\schurf_\lambda
	\cdot
\schurf_{\pas{\lambda_2,\lambda_3,\dots,\lambda_{m-1}}} = \\
\schurf_{\pas{\lambda_2,\lambda_3,\dots,\lambda_{m}}}
	\cdot
\schurf_{\pas{\lambda_1,\lambda_2,\dots,\lambda_{m-1}}} -
\schurf_{\pas{\lambda_2-1,\lambda_3-1,\dots,\lambda_{m}-1}}
	\cdot
\schurf_{\pas{\lambda_1+1,\lambda_2+1,\dots,\lambda_{m-1}+1}}
\end{multline}
\end{cor}
To provide a combinatorial proof for (a special case of) this identity was the
initial point for the work in \cite{fulmek:2001}.
\begin{proof}
Consider the matrix 
$\jtmat{\pas{\lambda+\constpart{1}{m}}/{\constpart{1}{m}}}$ (which is equal to
$\jtmat{\lambda}$). Then
\eqref{eq:Dodgson} translates to
\begin{multline*}
\schurf_{\pas{\lambda+\constpart{1}{m}}/\constpart{1}{m}}
	\cdot
\schurf_{\deletepart{\addtoparts{\lambda}{1}}{1,m}/\deletepart{\constpart{1}{m}}{1,m}}
	=\\
\schurf_{\deletepart{\addtoparts{\lambda}{1}}{1}/\deletepart{\constpart{1}{m}}{1}}
	\cdot
\schurf_{\deletepart{\addtoparts{\lambda}{1}}{m}/\deletepart{\constpart{1}{m}}{m}}
	-
\schurf_{\deletepart{\addtoparts{\lambda}{1}}{1}/\deletepart{\constpart{1}{m}}{m}}
	\cdot
\schurf_{\deletepart{\addtoparts{\lambda}{1}}{m}/\deletepart{\constpart{1}{m}}{1}}.
\end{multline*}
by \eqref{eq:JacobiTrudiSkew} and \eqref{eq:minorspartitions}. By
\eqref{eq:SchurShiftInvariance}
this simplifies to
$$
\schurf_\lambda\cdot\schurf_{\deletepart{\pas{\lambda+\constpart{1}{m}}}{1,m}}
	=
\schurf_{\deletepart{\pas{\lambda+\constpart{1}{m}}}{1}}\cdot\schurf_{\deletepart{\pas{\lambda}}{m}}
	-
\schurf_{\deletepart{\pas{\lambda}}{1}}\cdot\schurf_{\deletepart{\pas{\lambda+\constpart{1}{m}}}{m}},
$$
which is equivalent to \eqref{eq:DodgsonSchur}.
\end {proof}

So far, we proved the Schur function identity \eqref{eq:DodgsonSchur} by recognizing it
as a special case of the general determinantal identity \eqref{eq:Dodgson}. But this
works also the other way round:

\begin{obs}[Schur function identities imply equivalent determinantal identities]
\label{obs:genericSchurIdentity}
Assume that we have an identity $\identity{S}$ involving 
skew Schur functions.

Then by \eqref{eq:JacobiTrudiSkew} and Proposition~\ref{pro:hsym-alg-ind},
$\identity{S}$ translates to a 
determinantal identity $\identity{D}$ which is
\EM{equivalent} to $\identity{S}$ 
as follows:
\bit
\item Rewrite each skew Schur function in $\identity{S}$ as a determinant of complete homogeneous
	functions, according to \eqref{eq:JacobiTrudiSkew},
\item and then replace each entry $\csymf_{r}$ by variable 
$y_{r}$, where $\pas{y_r}_{r=0}^\infty$ is a set of independent variables,
according to Proposition~\ref{pro:hsym-alg-ind}.
\eit
\end{obs}

We call (skew) Schur function identities which are valid for \EM{arbitrary} shapes
$\lambda/\mu$ \EM{generic} (skew) Schur function identities: Note that
\eqref{eq:DodgsonSchur} is a generic identity in this sense. All identities we
shall consider in the rest of this paper are generic.

In particular, we may apply \eqref{eq:DodgsonSchur} to $\lambda=\pas{\lambda_1,\dots,\lambda_m}$,
where $\lambda_j=\pas{m-j+1}\cdot \pas{m}$. It is easy to see that for this choice of $\lambda$,
all the entries
in $\pas{\csymf_{\lambda_j-j+i}}_{i,j=1}^{m}$ are \EM{distinct}, whence we may
replace them by independent variables (by Proposition~\ref{pro:hsym-alg-ind}).
This means that \eqref{eq:DodgsonSchur} implies the \EM{general} determinantal
identity \eqref{eq:Dodgson}:
The identities are, in fact, \EM{equivalent}.

It is clear, that this phenomenon will apply also to other
determinantal identities: In this paper, we shall present \EM{generic}
Schur function identities which are equivalent to classical determinantal identities,
and give combinatorial proofs for these Schur function identities.

\section{Determinants as nonintersecting lattice paths}\label{sec:main}
The \lgt--method \cite{gessel-viennot:1998,lindstroem:1973,karlin:1988} is well--known:
But since we want to present a ``\nilp''--proof of the
Cauchy--Binet formula \eqref{eq:CauchyBinet} which involves a slight generalization, we repeat
this beautiful idea here in some detail, using the Jacobi--Trudi--determinant as an
illustrating example.

\subsection{Lattice paths}

Consider the square lattice $\mathbb Z^2$, i.e., the directed
graph with vertices $\mathbb Z\times\mathbb Z$ (we shall call them \EM{points}),
where the set of arcs consists of
\bit
\item horizontal arcs $\harc{\pas{j,k}}$ from $(j, k)$ to $(j+1, k)$ for $j, k\in\Z$,  and
\item vertical arcs $\varc{\pas{j,k}}$ from $(j, k)$ to $(j, k+1)$ for $j, k\in\Z$.
\eit
Assign to these arcs the following weights:
\begin{align*}
\weight\of{\varc{\pas{j,k}}} &\defeq 1\text{ (i.e., vertical arcs have weight $1$)}, \\
\weight\of{\harc{\pas{j,k}}} &\defeq x_k\text{ (i.e., horizontal arcs have weight $x_{\text{height of the arc}}$)}.
\end{align*}

A \EM{lattice path} $p$ of length $k$ connecting starting point $v$ to ending point $w$ is a sequence
of points $\pas{v=v_0,v_1,\dots,v_k=w}$, such that $\pas{v_{i-1},v_{i}}$ is a (horizontal
or  vertical) arc $a_i$ for $i=1,2,\dots,k$.
We say that all these arcs $a_i$ and points $v_i$ \EM{belong} to the path $p$ and write
$a_i\in p$ and $v_i\in p$. 

For some set $S$ of points (arcs) we denote by $p\cap S$ the set of points (arcs) in $S$
that belong to $p$. The weight of $p$ is defined as the product of the weights of the arcs belonging to $p$
\begin{equation}
\label{eq:pathweight}
\weight\of p=\prod_{a\in p}\weight\of a.
\end{equation}

Denote by $\paths{v}{w}$ the set of \EM{all} lattice paths connecting starting point
$v$ to ending point $w$.
Observe that we may view the complete homogeneous symmetric function
$$
\csymf_m\of{x_j,x_{j+1},\dots,x_k}, k\geq j
$$
as the \EM{generating function} of $\paths{v}{w}$, where
$v=\pas{t,j}$ and $w=\pas{t+m,k}$, with arbitrary $t\in\Z$:
$$
\gfpaths{v}{w}=\sum_{p}\weight\of p,
$$
where the sum is over all paths $p$ connecting $v$ and $w$.
(See \figref{fig:hsymgflp}.)

Note that the ending point of some path $p$
never lies below its
starting point. Later, we shall also consider \EM{trails} in the \EM{undirected} graph
$\mathbb Z^2$: A \EM{trail} $p^\prime$ of length $k$ connecting point $v$ to point $w$ is a sequence
of points $\pas{v=v_0,v_1,\dots,v_k=w}$, such that $\pas{v_{i-1},v_{i}}$ \EM{or}
$\pas{v_{i},v_{i-1}}$ is a (horizontal
or  vertical) arc $a_i$ for $i=1,2,\dots,k$ (i.e., $p^\prime$ may use arcs ``in the wrong direction'').
For such trail, it is
not clear whether $v$ or $w$ is the starting or ending point, so in order to
avoid confusion, for some lattice path $p$ we shall call
\bit
\item the starting point $p$ its \EM{lower} point,
\item and the ending point $p$ its \EM{upper} point
\eit
from now on.

\begin{figure}
\caption{Complete homogeneous symmetric functions may be viewed as generating functions of
lattice paths with fixed starting (lower) and ending (upper) point. The picture below illustrates this
for $\csymf_8\of{x_5,x_6,\dots,x_{13}}$, which appears as the sum of the weights of all
lattice paths connecting starting point $\pas{-4,5}$ to ending point $\pas{4,13}$, showing
the lattice paths associated to the two monomials
$\pas{x_5\cdot x_6^4\cdot  x_7\cdot  x_9^2}$ and
$\pas{x_{10}\cdot  x_{11}^2\cdot  x_{12}^4\cdot  x_{13}}$
in $\csymf_8\of{x_5,x_6,\dots,x_{13}}$. Note that we may shift the picture horizontally
by an arbitrary vector $\pas{t,0}$, $t\in\Z$, without changing the generating function
of the lattice paths.
}
\label{fig:hsymgflp}
\input graphics/hsymgflp
\end{figure}

\subsection{Young tableaux and \nilp s}

The \lgt\ interpretation gives an equivalent
description of a semistandard Young tableau $\tableau$ of shape $\lambda/\mu$ as an
$m$--tuple $P=\pas{p_1,\dots,p_{m}}$ of \nilp s, where $m\defeq\length{\lambda/\mu}$, as follows.
The $i$--th path $p_i$ starts at lower point $(\mu_i-i, 1)$ and ends at upper point
$(\lambda_{i}-i, \n)$: We call these points the lower/upper points associated to the
skew shape $\lambda/\mu$, 
and we count these points always \EM{from the right}.
The $k$--th \EM{horizontal} step in $p_i$ goes from
$\pas{\mu_i-i+k-1,x}$ to $\pas{\mu_i-i+k,x}$, where $x$ is the $k$--th entry in row $i$ of
$\tableau$. 

\begin{figure}
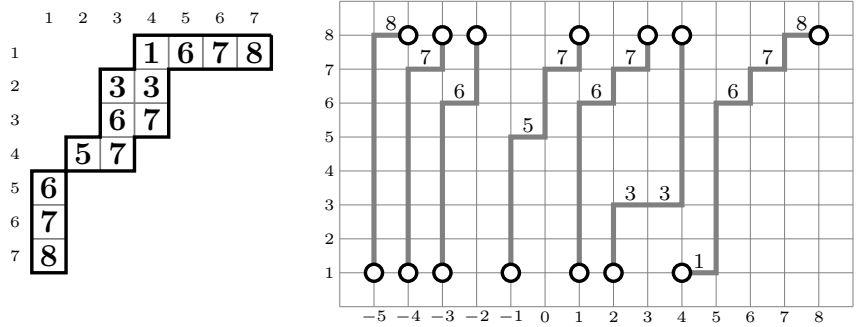

\caption{The left picture presents a semistandard Young tableau $\tableau$ of shape $\lambda/\mu$, where
$\lambda=\pas{9, 6, 6, 5, 3, 3, 3}$ and $\mu=\pas{5, 4, 4, 3, 2, 2, 2}$. Assuming that the
entries of $\tableau$ 
are chosen from $\setof{1,2,\dots,8}$ (i.e.: $\tableau$ is an $8$--semistandard Young tableau), the right picture
shows the corresponding family of $7=\length{\lambda/\mu}$ nonintersecting lattice paths: Note that the height
of the $j$--th horizontal step in the $i$--th path (the paths are counted from right to left) is equal
to the $j$--th entry in row $i$ of $\tableau$.}
\label{fig:skew-Young-tableau}
\input graphics/ferrers
\end{figure}

Note that the conditions on the entries of $\tableau$ imply that
no two paths $p_i$ and $p_j$ thus defined have a lattice point in common: Such an
$m$-tuple of paths is called \EM{nonintersecting},
see the right picture of Figure~\ref{fig:skew-Young-tableau} for an illustration.
An $m$--tuple of paths which is \EM{not} nonintersecting is called \EM{intersecting}.

In fact, this translation of tableaux to nonintersecting lattice paths is
a \EM{bijection} between the set of all $\n$--semistandard Young tableaux of shape $\lambda/\mu$
and the set of all $\length{\lambda/\mu}$--tuples of nonintersecting lattice paths with lower (starting) and
upper (ending) points as defined above.

This bijection is \EM{weight preserving} if we define
the weight of an $m$--tuple $P=\pas{p_1,\dots,p_m}$ of lattice paths in the obvious way, i.e., as
\begin{equation}
\label{eq:dfnweightlp}
\weight\of P \defeq\prod_{k = 1}^{m} \weight\of{p_k} = \prod_{k = 1}^{\n} x_k^{\nofentries{P}{k}},
\end{equation}
where $\nofentries{P}{k}$ is the number of horizontal steps at height $k$ in $P$. So in definition
\eqref{eq:skewSchur} we could equivalently replace symbol ``$\tableau$'' by symbol ``$P$'',
and sum over $m$--tuples of nonintersecting lattice
paths with prescribed lower and upper points instead of tableaux with prescribed shape.

Note that the horizontal coordinates of lower and upper points determine uniquely the
\EM{shape} $\lambda/\mu$ of the tableau, and the vertical coordinate (we shall call the 
vertical coordinate of points the \EM{level} in the following) of the ending points determines
uniquely the \EM{set of entries} $\setof{1,2,\dots,\n}$ of the tableau. (The choice of the
shift parameter $t$ does influence neither the shape nor the set of entries.)

Observe that the operation of deleting part $i$ in $\mu=\pas{\mu_1,\dots,\mu_m}$ and part $j$ in
$\lambda=\pas{\lambda_1,\dots,\lambda_m}$ (defined in \eqref{eq:deletepart}) translates
to the removal of the $i$--th lower point and the $j$--th upper point associated to
$\lambda/\mu$,
and by \eqref{eq:minorspartitions}, this removal of lower/upper points translates to deleting row $i$ and column $j$ in $\jtmat{\lambda/\mu}$. Clearly,
this generalizes to the following observation:
\begin{obs}
\label{obs:deletepointminors}
Minors of $\jtmat{\lambda/\mu}$
consisting of rows $\setof{i_1,\dots,i_k}$ and columns $\setof{j_1,\dots,j_l}$
are in one--to--one correspondence to the selection of
\bit
\item lower points with indices in $\setof{i_1,\dots,i_k}$
\item and upper points with indices in $\setof{j_1,\dots,j_l}$
\eit
from the points associated to the shape $\lambda/\mu$.
\end{obs}

\subsection{The \lgt--proof of the Jacobi--Trudi identity}
Denote the $i$--the lower point $(\mu_i-i, 1)$ by $s_i$, and the $j$--th ending point
$(\lambda_{j}-j, \n)$ by $t_j$.
Observe that the entry $\pas{i,j}$ in the matrix $\jtmat{\lambda/\mu}$ is the
\EM{generating function} 
$$
\gfpaths{s_i}{t_j}\defeq\sum_{p}\weight\of{p}=\csymf_{\lambda_j-\mu_i-j+i},
$$
where the sum is over all lattice paths $p$ that run from $s_i$
to $t_j$,
and where the weight $\weight\of{p}$ is defined as in \eqref{eq:pathweight}.
By expanding the determinant in
\eqref{eq:JacobiTrudiSkew}, we thus obtain
\begin{equation}
\label{eq:detgf}
\detof{\csymf_{\lambda_j-\mu_i-j+i}}_{i,j=1}^{m}
=
\sum_{\pi\in\symm_m}\sgn\of{\pi}\cdot
\prod_{j=1}^m
\pas{\gfpaths{s_{\pi_j}}{t_j}},
\end{equation}
where $m=\length{\lambda/\mu}$. 
Consider the following set
\begin{equation*}
\alljt_{\lambda/\mu}
\defeq
\bigcup_{\pi\in\symm_m}
\paths{s_{\pi_1}}{t_1}\times\paths{s_{\pi_2}}{t_2}\times\cdots\times\paths{s_{\pi_m}}{t_m}
\end{equation*}
of $m$--tuples $P$ of lattice paths connecting the permuted lower points
$\pas{s_{\pi_1},s_{\pi_2},\dots,s_{\pi_m}}$
with the
upper points
$\pas{t_1,t_2,\dots,t_m}$,
where in addition to the weight $\weight\of P$, the elements of $\alljt_{\lambda/\mu}$
also carry a \EM{sign} which equals the sign of the respective permutation $\pi$:
$$\sgn\of P = \sgn\of\pi.$$
Then \eqref{eq:detgf} can be rewritten equivalently as the generating function of 
$\alljt_{\lambda/\mu}$:
\begin{equation}
\label{eq:detgf2}
\detof{\csymf_{\lambda_j-\mu_i-j+i}}_{i,j=1}^{m}
=
\sum_{P\in{\alljt_{\lambda/\mu}}}\sgn\of P\cdot \weight\of P.
\end{equation}
Denote by $\nijt_{\lambda/\mu}$ the subset of \EM{nonintersecting} $m$-tuples of lattice
paths in $\alljt_{\lambda/\mu}$. In order to prove \eqref{eq:JacobiTrudiSkew}, we only
need to show that we do in fact have
\begin{equation}
\label{eq:detisgf}
\detof{\csymf_{\lambda_j-\mu_i-j+i}}_{i,j=1}^{m}
=
\sum_{P\in{\nijt_{\lambda/\mu}}}\sgn\of P\cdot \weight\of P, 
\end{equation}
since $P\in{\nijt_{\lambda/\mu}}$ implies $\pi=\idperm$ (i.e., all these ``surviving''
objects have sign $+1$) and $\schurf_{\lambda/\mu}=\sum_{P\in{\nijt_{\lambda/\mu}}}\sgn\of P\cdot \weight\of P$
by definition \eqref{eq:skewSchur}.
This certainly can be achieved by showing that all the signed weights of $m$--tuples
$P$ in $\ijt_{\lambda/\mu}\defeq\alljt_{\lambda/\mu}\setminus\nijt_{\lambda/\mu}$
\EM{cancel} in \eqref{eq:detgf2}. To this end,
we shall present an
\EM{involution} (a self--inverse bijective mapping) on \EM{intersecting} $m$--tuples
of lattice paths
$$
\invol:
\pas{\ijt_{\lambda/\mu}}
\to
\pas{\ijt_{\lambda/\mu}}
$$
which is
\bit
\item \EM{weight--preserving}, i.e., $\weight\of{\invol\of{P}} = \weight\of P$
\item and \EM{sign--reversing}, i.e., $\sgn\of{\invol\of{P}} = -\sgn\of P$.
\eit
The construction of $\invol$ is simple: Let $P$ be an \EM{intersecting} $m$--tuple in
$\ijt_{\lambda/\mu}$. Consider the \EM{smallest} point of
intersection $q$ in $P$, in lexicographic order:
$$
\pas{a,b}\leq\pas{c,d}:\Leftrightarrow a\leq b \vee
\pas{a=b\wedge b\leq c}.
$$
Observe that there are precisely \EM{two} paths $p_k$ and $p_l$ meeting in $q$
(by the minimality of $q$): $\invol\of P$ is obtained from $P$ by \EM{interchanging}
the initial segments (from lower points up to $q$) of $p_k$ and $p_l$, see
\figref{fig:lgt-involution} for an illustration.

It is immediately clear that $\invol$
is an involution which is weight--preserving and sign--reversing (since it modifies the
original permutation associated to $P$ by the transposition corresponding to the swapping
of the lower points of $p_k$ and $p_l$).

Clearly, $\invol$
describes the pairwise cancellation
$$\sgn\of{P}\cdot\weight\of{P}+\invol\of{\sgn\of{P}\cdot\weight\of{P}} = 0$$
of all objects from $\ijt_{\lambda/\mu}$ in the
sum \eqref{eq:detgf2}.
Stated otherwise: Only the \EM{nonintersecting} objects in $\nijt_{\lambda/\mu}\subseteq\alljt_{\lambda/\mu}$
``survive'', which proves \eqref{eq:detisgf}.
\hfill\qedsymbol

\begin{figure}
\caption{Illustration of the \lgt--involution $\invol$. Both pictures show quadruples of
lattice paths belonging to $\ijt_{\lambda/\mu}$, where $\lambda=\pas{8,6,4,3}$ and
$\mu=\pas0$. There are $4$ points of intersections at positions
$$
\pas{6,2} < \pas{6,5} < \pas{9,5} < \pas{9,6},
$$
written in lexicographic order. 
(Note that the lattice paths are drawn with rounded corners
and small offsets here, just to make obvious the run of the paths, which is not clear for
intersecting paths). The smallest such point in lexicographic order is
$q\defeq\pas{6,2}$ (indicated by a circle). The right picture is obtained from the left picture by interchanging
the initial segments (from lower points up to $q$) of the two paths intersecting in $q$, and
vice versa. To the left picture, the identity permutation (i.e., sign $+1$) is associated, while to the right picture, the transposition $\pas{3,4}$ (i.e., sign $-1$) is associated. This illustrates that in the determinantal
expansion of $\schurf_\lambda=\detof{\csymf_{\lambda_j-j+i}}$, the following terms cancel:
$$
\pas{x_7^3}\pas{x_5^4}\pas{x_2^3 x_5^3}\pas{x_1^2 x_2 x_3^2 x_6^3}
-
\pas{x_7^3}\pas{x_5^4}\pas{x_2^4 x_3^2 x_6^3}\pas{x_1^2 x_5^3}
=0.
$$
}
\label{fig:lgt-involution}
\input graphics/lgt-involution
\end{figure}

\subsection{Viewing determinants as nonintersecting lattice paths} Our considerations
so far showed that
the \EM{generic} determinant $\pas{y_{i,j}}_{i,j=1}^m$ (here, generic means that the entries are
independent variables) may be viewed as a skew Schur function of \EM{appropriate} shape
$\lambda/\mu$ (here, appropriate means all entries in $\jtmat{\lambda/\mu}$ are distinct),
hence as the
\EM{generating function of \nilp s}.

We shall demonstrate the power of this point of view by giving a simple proof of the
Cauchy--Binet formula~\eqref{eq:CauchyBinet}, which will become even more transparent
if we prove the multiplicativity of the determinant function as a preparatory step.


\subsubsection{Multiplicativity of the determinant function: A proof ``by example''.}
Given some arc $a$ in the lattice $\Z^2$, we say that a path $p$ \EM{starts} in $a$
(has $a$ as its \EM{lower} arc), if $p$ starts in the upper (if $a$ is vertical) or
right (if $a$ is horizontal) point of $a$. Likewise, we say that $p$ \EM{ends} in
$a$ (has $a$ as its \EM{upper} arc), if $p$ ends in the lower (if $a$ is vertical) or
left (if $a$ is horizontal) point of $a$.

We want to illustrate the multiplicativity of the determinant function
\begin{equation}
\label{eq:detmultfunc}
\detof{a\cdot b} = \detof{a}\cdot\detof{b}
\end{equation}
by considering the special case
$$
a\defeq\pas{\csymf_{\lambda_j-j+i}\of{x_1,x_2,\dots x_7}}_{i,j=1}^4\text{ and }
b\defeq\pas{\csymf_{\sigma_j-\lambda_i-j+i}\of{x_8,x_9\dots x_{13}}}_{i,j=1}^4,
$$
where $\lambda=\pas{8,6,4,3}$ and $\sigma=\pas{18,16,13,11}$ .

By the \lgt--interpretation, we may view
$\detof{a}$ as the generating function of quadruples of nonintersecting
lattice paths connecting lower points $\asvector{r}=\pas{r_1,r_2,r_3,r_4}$ with upper arcs
$\asvector{s}=\pas{s_1,s_2,s_3,s_4}$, where
\begin{align*}
\asvector{r} &= \pas{{
	\pas{-1,1},\;\pas{-2,1},\;\pas{-3,1},\;\pas{-4,1}}},\\
\asvector{s} &= \pas{{
	\pas{\pas{7,7},\pas{7,8}},\;\pas{\pas{4,7},\pas{4,8}},\;\pas{\pas{1,7},\pas{1,8}},\;\pas{\pas{-1,7},\pas{-1,8}}}}.
\end{align*}

Likewise, we may view $\detof{b}$ as the generating function of quadruples of nonintersecting
lattice paths connecting lower arcs $\asvector{s}$ (as above: Note that the set of variables involved in $b$ is $\setof{x_8,\dots,x_{13}}$!) with ending
points $\asvector{t}=\pas{t_1,t_2,t_3,t_4}$, where
$$
\asvector{t} = \pas{{
	\pas{17,13},\;\pas{14,13},\;\pas{10,13},\;\pas{7,13}}}.
$$
See \figref{fig:DeterminantMultiplicative}
for an illustration.

\begin{figure}
\caption{Illustration: Multiplicativity of the determinant.
The picture shows the
lattice paths corresponding to the monomial
$$
\scriptstyle
\pas{\pas{x_7^3}\cdot\pas{x_9^2 x_{11} x_{13}^5}}
\cdot
\pas{\pas{x_5^4}\cdot\pas{x_7^2 x_{11}^5 x_{12}^2}}
\cdot
\pas{\pas{x_2^3 x_4^3}\cdot\pas{x_9^5 x_{11}^5}}
\cdot
\pas{\pas{x_1^3 x_3^4 x_5}\cdot\pas{x_8^4 x_{10}^6}},\hskip5em
$$
which appears in the expansion of the product of skew Schur functions
$$\schurf_\lambda\of{x_1,x_2,\dots x_7}\cdot
\schurf_{\sigma/\lambda}\of{x_8,x_9,\dots x_{13}},$$
where $\lambda=\pas{8,6,4,3}$ and $\sigma=\pas{18,16,13,11}$.
}
\label{fig:DeterminantMultiplicative}
\input graphics/DeterminantMultiplicative
\end{figure}

So the $\pas{i,j}$--entry in $a\cdot b$ is
$$
\sum_{k=1}^4\gfpaths{r_i}{s_k}\cdot\gfpaths{s_k}{t_j},
$$
which may be viewed as the generating function $\gfof{\pathsmod{r_i}{t_j}{\prime}}$ of
the following set of \EM{constrained} paths (see \figref{fig:lgt-modified} for an illustration):
$$
\pathsmod{r_i}{t_j}{\prime}\defeq
\bigcup_{k=1}^4 \setof{p:\text{ $p$ is a lattice path connecting $r_i$ to $t_j$ \EM{passing through} $s_k$}}.
$$

\begin{figure}
\caption{Illustration: The \lgt--involution for the ``constrained'' lattice paths
$\gfof{\pathsmod{r_i}{t_j}{\prime}}$: The fact that \EM{every} constrained path must pass through
one of the arcs $s_1,s_2,s_3$ or $s_4$ is indicated 
by omitting
all \EM{other} arcs connecting level $7$ and level $8$.
The picture shows a path connecting $r_1$ to $t_1$ which passes through the arc $s_2$,
and a path connecting $r_4$ to $t_3$,
which passes through the arc $s_3$. The smallest point of intersection in lexicographic order is
$q\defeq\pas{-4,9}$ (indicated by a circle;
as in \figref{fig:lgt-involution}, the lattice paths are drawn with rounded corners
here, just to make obvious the run of the paths.) Clearly, the \lgt--involution gives
a path connecting $r_1$ to $t_3$ which passes through $s_2$, and a path connecting $r_4$ to $t_1$
which passes through $s_3$.
}
\label{fig:lgt-modified}
\input graphics/lgt-modified
\end{figure}

So as in the above proof of the Jacobi--Trudi identity, $\detof{a\cdot b}$ can be rewritten
equivalently as the generating function of the set ${\mathfrak D}^\prime$ of quadrupels of
\EM{constrained} paths,
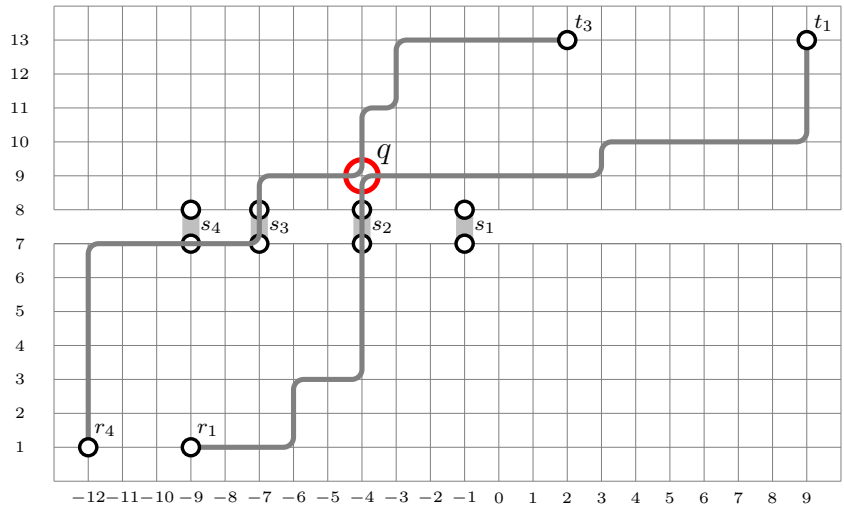
\begin{equation*}
{\mathfrak D}^\prime\defeq
\bigcup_{\pi\in\symm_4}
\pathsmod{s_{\pi_1}}{t_1}{\prime}\times
\pathsmod{s_{\pi_2}}{t_2}{\prime}\times
\pathsmod{s_{\pi_3}}{t_3}{\prime}\times
\pathsmod{s_{\pi_4}}{t_4}{\prime},
\end{equation*}
i.e., as
\begin{equation}
\label{eq:detmultgf}
\detof{a\cdot b}
=
\sum_{P\in{\alljt^\prime}}\sgn\of P\cdot \weight\of P.
\end{equation}
Denote by $\nijt^\prime$ the subset of \EM{nonintersecting} $m$-tuples of constrained lattice
paths in $\alljt^\prime$. Clearly, the \lgt--operation $\invol$ can be
applied 
to \EM{intersecting} $m$-tuples of constrained lattice paths (see \figref{fig:lgt-modified} for an illustration)
and establishes a weight--preserving and sign--reversing
involution, i.e.,
$$\detof{a\cdot b}=\gfof{\nijt^\prime}.$$
On the other hand, $\nijt^\prime$ appears (just \EM{look} at \figref{fig:DeterminantMultiplicative}!)
as the Cartesian product of
\bit
\item the set
of nonintersecting quadrupels of lattice paths connecting $\asvector{r}$ and $\asvector{s}$
\item and the set of nonintersecting quadrupels of lattice paths connecting $\asvector{s}$ and $\asvector{t}$,
\eit
i.e.,
$$
\gfof{\nijt^\prime}=\detof a\cdot \detof b.
$$

We presented our argument for a specific choice of partitions $\lambda\subpart\sigma$
to make it more tangible, but it is clear that it holds for arbitrary choices of (semi)partitions.
%
%
%
%
So fix $m>0$ and consider the following partitions $\lambda$ and $\sigma$,
$\length{\lambda}=\length{\sigma}=m$:
\begin{align*}
\lambda&=\pas{(m\cdot\pas{m-1},\pas{m-1}\cdot\pas{m-1},\dots,\pas{m-1}} 
\text{ (i.e., } \lambda_j=\pas{m-j+1}\cdot\pas{m-1}\text{),}\\
\sigma &=\constpart{\pas{m\cdot m}}{m}=\pas{m^2,m^2,\dots,m^2}
\text{ (i.e., } \sigma_i=m^2\text{).}
\end{align*}
Then the entries in
$a\defeq\jtmat{\lambda}\of{x_1,\dots,x_\n}$ and in $b\defeq\jtmat{\sigma/\lambda}\of{x_{\n+1},\dots,x_{2\n}}$ are all distinct. Hence (recall Observation~\ref{obs:genericSchurIdentity})
the identity
$$
\detof{a\cdot b} = \detof a \cdot \detof b
$$
implies the mulitplicativity of the determinant function \EM{in general}.
\hfill\qedsymbol

\subsubsection{The Cauchy--Binet formula: Another proof ``by example''.} Now it is ``almost
immediate'' to \EM{see} the Cauchy--Binet formula as an obvious generalization of the multiplicativity of the determinant.

\begin{thm}[Cauchy--Binet formula]
Let $a$ be an $m\times n$--matrix and $b$ be an $n \times m$--matrix. Then we have
\begin{equation}
\label{eq:CauchyBinet}
\detof{a\cdot b} = \sum_{S\subseteq\brk{n},\cardof{S}=m}
	\detof{\minor{a}{\brk{m}}{S}}\cdot\detof{\minor{b}{S}{\brk{m}}}.
\end{equation}
\end{thm}
Note that this formula holds trivially if $m>n$ (since the determinant
$\detof{a\cdot b}$ is zero in this case, as is the empty sum in \eqref{eq:CauchyBinet}), and
amounts to the multiplicativity of the determinant if $m=n$.

We illustrate this formula by a special case and consider the matrices
$$
a\defeq\pas{\csymf_{\lambda_j-j+i}\of{x_1,\dots x_7}}_{\pas{i,j}\in\pas{\brk{3}\times\brk{4}}}\text{ and }
b\defeq\pas{\csymf_{\sigma_j-\lambda_i-j+i}\of{x_8,\dots x_{13}}}_{\pas{i,j}\in\pas{\brk{4}\times\brk{3}}},
$$
where $\lambda=\pas{9,7,5,4}$ and $\sigma=\pas{19,17,14}$.

Again, we want to employ the \lgt--interpretation.
To this end, we consider the vectors of points $\asvector{r}=\pas{r_1,r_2,r_3}$ and
$\asvector{t}=\pas{t_1,t_2,t_3}$, and the vector of arcs,
$\asvector{s}=\pas{s_1,s_2,s_3,s_4}$, where
\begin{align*}
\asvector{r} &= \pas{{
	\pas{-1,1},\;\pas{-2,1},\;\pas{-3,1}}},\\
\asvector{s} &= \pas{{
	\pas{\pas{8,7},\pas{8,8}},\;\pas{\pas{5,7},\pas{5,8}},\;
	\pas{\pas{2,7},\pas{2,8}},\;\pas{\pas{0,7},\pas{0,8}}}},\\
\asvector{t} &= \pas{{
	\pas{18,13},\;\pas{15,13},\;\pas{11,13}}},
\end{align*}
(See \figref{fig:Cauchy-Binet} for an illustration.)

\begin{figure}
\caption{Illustration: The Cauchy--Binet formula. Triples of nonintersecting
lattice paths connecting
points $\asvector{r}=\pas{r_1,r_2,r_3}$ and $\asvector{t}=\pas{t_1,t_2,t_3}$, where
each single path must pass through an arc in $\asvector{s}=\pas{s_1,s_2,s_3,s_4}$,
can be ``cut in two halves'' along $\asvector{s}$. Note that the arcs used by the paths
constitute a $3$--element subvector $\asvector{s}^\prime$, and the
halves appear as tripels of nonintersecting
lattice paths $\nijt\of{\asvector{r}\to\asvector{s}^\prime}$ and
$\nijt\of{\asvector{s}^\prime\to\asvector{t}}$,
respectively. In the picture, we have $\asvector{s}^\prime=
\pas{s_1,s_2,s_4}$.
}
\label{fig:Cauchy-Binet}
\input graphics/Cauchy-Binet
\end{figure}

As before,
observe that the
$\pas{i,j}$--entry in $a\cdot b$
may be viewed as the generating function $\gfpaths{r_i}{t_j}^\prime$ of
\EM{constrained} paths which must pass through one arc in $\asvector{s}$. By the \lgt--involution, the determinant $\detof{a\cdot b}$
appears as the generating function of quadrupels of nonintersecting labeled paths
$\nijt^\prime\of{\asvector{r}\to\asvector{t}}$, which (as above --- just \EM{look} at \figref{fig:Cauchy-Binet}!) appears as
$$
\nijt^\prime\of{\asvector{r}\to\asvector{t}} = \bigcup_{\asvector{s}^\prime}
	\pas{
		\nijt\of{\asvector{r}\to\asvector{s}^\prime}
		\times
		\nijt\of{\asvector{s}^\prime\to\asvector{t}},
	}
$$
where the union runs over all $3$--element subvectors $\asvector{s}^\prime$ of
$\asvector{s}$. By the same reasoning as above, this shows \eqref{eq:CauchyBinet}.
\hfill\qedsymbol
%
%
%

\section{Products of determinants as overlays of lattice paths}\label{sec:overlays}
In the following, all skew Schur functions are considered as functions of the same set
of variables
$\pas{x_1,\dots,x_{\n}}$. (Equivalently, all tableaux have entries from the set
$\setof{1,\dots,\n}$, and all families of \nilp s have lower
points on level $1$ and upper
points on level $\n$).

By \eqref{eq:skewSchur} and the \lgt--interpretation of Young tableaux as \nilp s,
we may view the product of two skew Schur functions as the generating function of
pairs of $m$--tuples of \nilp s
$$
\schurf_{\lambda/\mu}\cdot\schurf_{\sigma/\tau} =
	\sum_{\pas{P_1,P_2}}\weight\of{P_1}\cdot\weight\of{P_2}.
$$
Here, the sum runs over all pairs $\pas{P_1,P_2}$, where $P_1$ is a $\length{\lambda/\mu}$--tuple
of \nilp s corresponding to the shape $\lambda/\mu$ and $P_2$ is a
$\length{\sigma/\tau}$--tuple
of \nilp s corresponding to the shape $\sigma/\tau$.
Imagine that the lattice paths of $P_1$ and $P_2$ are coloured red and green, respectively:
This will give an \EM{overlay} of lattice paths,
see \figref{fig:nilpOverlay} for an illustration.

\begin{figure}
\caption{Illustration: A green $7$--tuple of \nilp s, corresponding to
shape $\lambda/\mu$, and a red $7$--tuple of \nilp s, corresponding to
shape $\sigma/\tau$, constitute an \EM{overlay} of \nilp s. Here, $\sigma=\pas{7, 4, 4, 3, 1, 1, 1}$, $\tau=\pas{3, 2, 2, 1}$,
$\lambda=\sigma+\constpart{1}{7}$ and $\mu=\tau+\constpart{1}{7}$.
(Note that $\schurf_{\lambda/\mu} = \schurf_{\sigma/\tau}$.) The red and green paths are drawn with a slight offset for graphical
reasons, the colour red is indicated by dashed lines.
The picture also shows the Young tableaux corresponding to the red and green
$7$--tuples of \nilp s.
}
\label{fig:nilpOverlay}
\input graphics/nilpOverlay-slide
\end{figure}

Such overlays 
give rise to a bijective construction,
which (to the best of our knowledge) was first used by Goulden \cite{goulden:1988}.
The same construction was used
in \cite{fulmek:2001} to describe and prove a class of Schur function identities, special
cases of which imply Dodgson's condensation formula and the Pl\"ucker relations.


We shall present this construction by way of an example: Consider skew shapes
$\lambda/\mu$ and $\sigma/\tau$, where $\sigma=\pas{7, 4, 4, 3, 1, 1, 1}$, $\tau=\pas{3, 2, 2, 1}$,
$\lambda=\sigma+\constpart{1}{7}$ and $\mu=\tau+\constpart{1}{7}$,
see \figref{fig:nilpOverlay} for an illustration.


Now consider the arcs, lower points and upper points of some pair $\pas{P_1,P_2}$ of 
$7$--tuples of lattice paths,
where $P_1$ corresponds to shape $\lambda/\mu$ and $P_2$ corresponds to
shape  $\sigma/\tau$.
Note that there may be arcs/points coloured both green \EM{and} red. Such arcs/points will
never be affected by the following constructions: We call them \EM{uncoloured}; the \EM{remaining} arcs/points
(which are coloured \EM{either} red \EM{or}
green) are called the \EM{coloured} arcs/points.

%
%
We construct \EM{bicoloured trails}
\bit
\item connecting coloured points
\item and using (only) coloured arcs
\eit
by the following algorithm: 

\begingroup
\leftskip1em
\rightskip1em
We start at some coloured point $s$ and identify the
\EM{unique} coloured arc $a$ incident with $s$ which is \EM{of the same colour} as $s$.
Then we follow the lattice path starting in $a$ in the implied direction (i.e., either
up/right if
$s$ is a lower point, or down/left if $s$ is an upper point).

Whenever we meet \EM{another} path on our way (necessarily, this path is of the other colour),
we ``change colour and direction'', i.e., we follow this new path \EM{and} change the direction (i.e., if we were moving up/right along the old path, we move down/left
along the new path, and vice versa). Note that such change of colour and orientation might
also occur at the very beginning: For instance, if $s$ is a green upper point, but there
are two red arcs incident with $s$, then we follow the red arc in the up/right direction. 

We stop if there is no possibility to go further,
i.e., if we end in another coloured point.

\endgroup



\figref{fig:bicolouredPaths-slide} illustrates this construction (see also \cite{fulmek:2001}).

\begin{figure}
\caption{
The picture shows the bicoloured trails (as thick grey trails with rounded corners)
for the example from \figref{fig:nilpOverlay}.}
\label{fig:bicolouredPaths-slide}
\input graphics/bicolouredPaths-slide
\end{figure}

The following observations are immediate from the construction:
\begin{obs}[{\bf Bicoloured trails always exist}]
\label{obs:arbitrary-point}
For \EM{every} coloured point $s$, there \EM{exists} a bicoloured trail starting at $s$.
\end{obs}

\begin{figure}
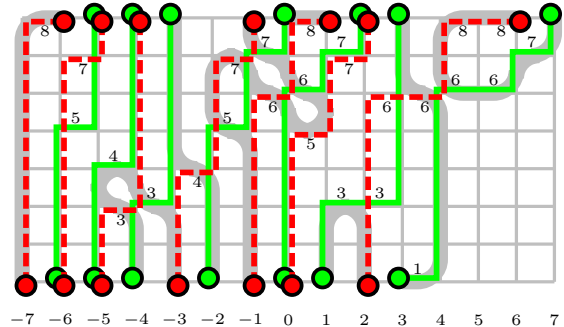

\caption{Illustration: Bicoloured trails can never cross. Since bicoloured
trails never use arcs coloured both red \EM{and} green, a meeting point of bicoloured trails
can (up to interchanging colours red and green) only occur in the situations shown below. By 
construction, the bicoloured trails necessarily run as indicated by the grey lines.}
\label{fig:noncrossing}
\input graphics/noncrossing
\end{figure}

\begin{obs}[{\bf Bicoloured trail can never cross}]
\label{obs:different-parity}
Bicoloured trails may have lattice points in common (they may intersect),
but they can never \EM{cross} (see \figref{fig:noncrossing}).
\end{obs}

Now consider some bicoloured trail $b$ in the overlay of \nilp s $\pas{P_1,P_2}$: 
Changing colours (green to red and vice versa)
\bit
\item of both ending points of $b$
\item and of all arcs of $b$
\eit
gives a new overlay of \nilp s  $\pas{P_1^\prime,P_2^\prime}$
(with different lower/upper points, see \figref{fig:recolouringPaths-slide}).
Clearly, we have:
\begin{obs}[{\bf Recolouring bicoloured trails is a weight preserving involution}]
\label{obs:weight-preserving-involution}
The recolouring of a bicoloured trail $b$ in an overlay  of \nilp s $\pas{P_1,P_2}$ is an \EM{involutive operation}, i.e., if we obtain
the overlay $\pas{P_1^\prime,P_2^\prime}$ by recolouring $b$ in $\pas{P_1,P_2}$, then
recolouring $b$ \EM{again} in $\pas{P_1^\prime,P_2^\prime}$ yields the original
$\pas{P_1,P_2}$. Moreover, this operation \EM{preserves 
the respective weights}, i.e.,
$$
\weight\of{P_1}\cdot\weight\of{P_2}=\weight\of{P_1^\prime}\cdot\weight\of{P_2^\prime}.
$$
\end{obs}
\begin{figure}
\caption{Illustration. Recolouring two bicoloured trails (indicated by thick grey
lines) from \figref{fig:bicolouredPaths-slide} changes the colour of lower and/or
upper points, thus changing the corresponding shape (shown below the paths).}
\label{fig:recolouringPaths-slide}
\input graphics/recolouringPaths-slide
\end{figure}

\begin{figure}
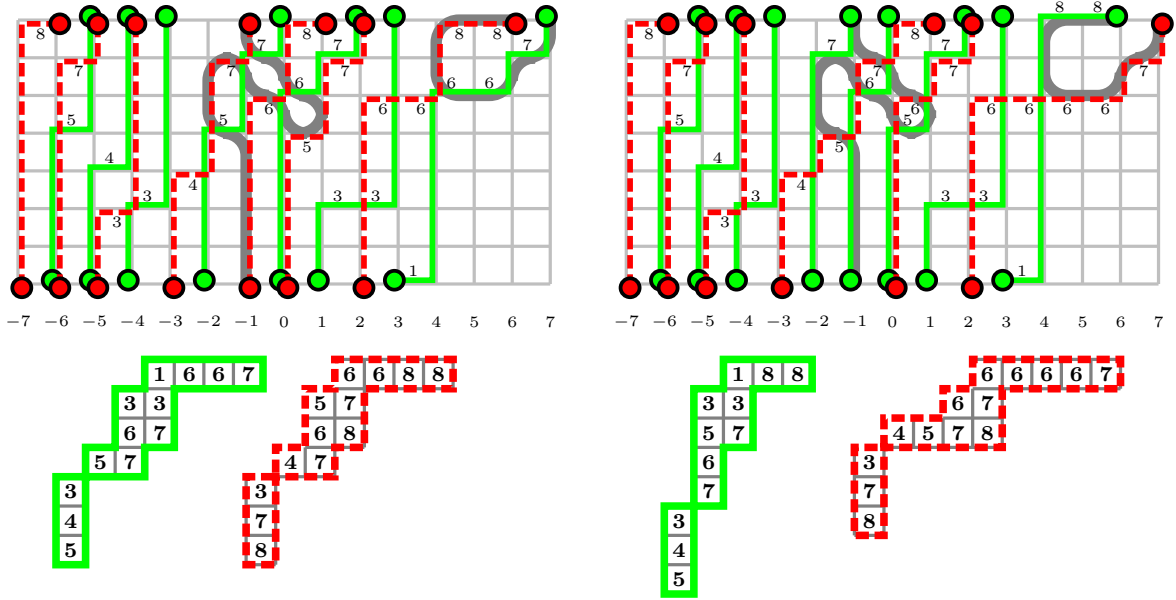

\caption{The ``radial orientation'' (indicated in the picture by triangles) of
coloured points depends on the position (upper or lower) and colour (red or green)
of the points.
}
\label{fig:bicolouredPaths2-slide}
\input graphics/bicolouredPaths2-slide
\end{figure}

Note that the operation of recolouring bicoloured trails changes the colours (red/green)
of coloured lower and/or upper points (which implies a change of the corresponding shapes,
see \figref{fig:recolouringPaths-slide}).
We want to encode this change in a convenient way: Imagine that all
lower/upper points are
arranged on a circle (see \figref{fig:bicolouredPaths2-slide}). Assign to coloured point $s$
the \EM{radial orientation} (with respect to this circle)
\bit
\item \EM{inwards}, if $s$ is a \EM{red upper point} or a \EM{green lower point},
\item \EM{outwards}, if $s$ is a \EM{green upper point} or a \EM{red lower point}.
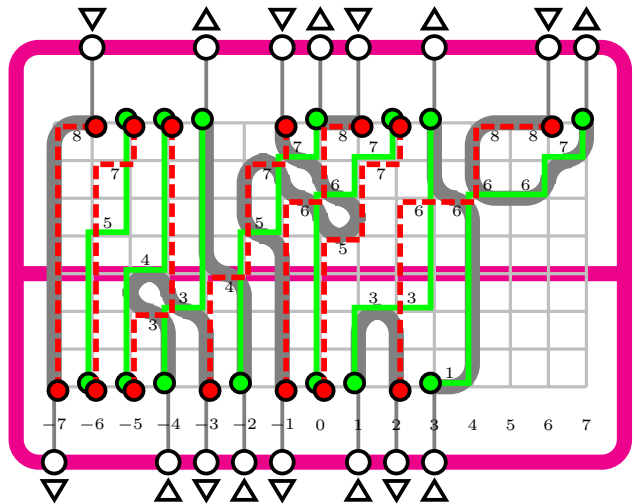
\eit


See \figref{fig:bicolouredPaths2-slide}. From the construction of bicoloured trails,
the following is immediate:
\begin{obs}[{\bf Bicoloured trails connect points of different radial orientation}]
\label{obs:different-orientation}
Bicoloured trails never connect points of the same radial orientation
(i.e., two points oriented both inwards or both outwards).
\end{obs}

\begin{figure}
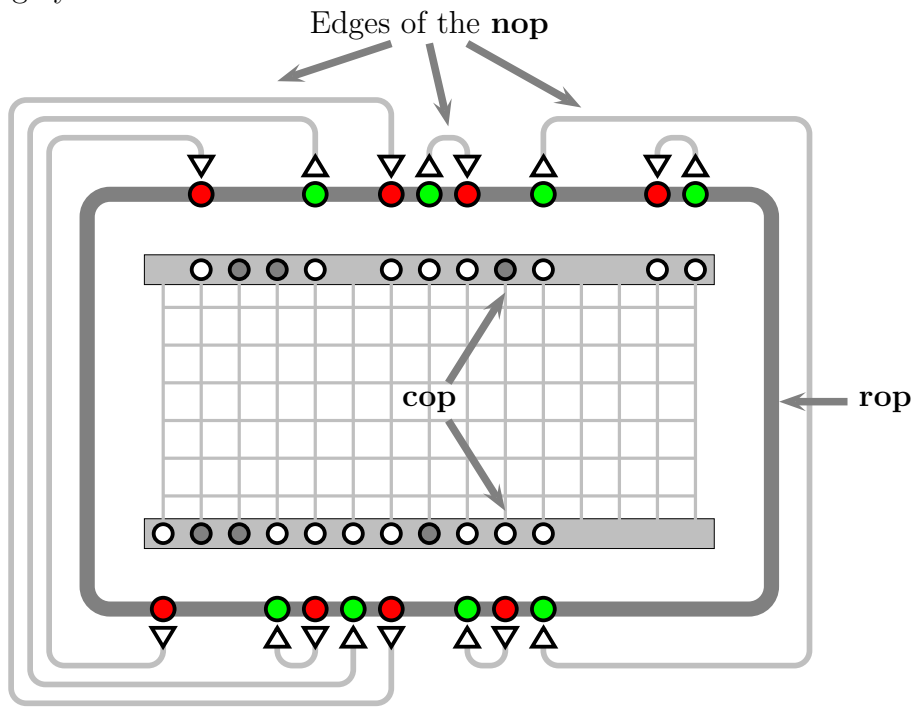

\caption{Coloured points may be green or red: The pattern of radial
orientations \pointarr\ (indicated by inward/outward pointing triangles)
encodes the actual colours of the coloured points in the \pointconf\ 
(indicated by white circles, the uncoloured points appear as grey circles)
and thus determines the corresponding shapes $\lambda/\mu$ and $\sigma/\tau$.
The picture illustrates this for the example from
\figref{fig:recolouringPaths-slide}: The bicoloured trails constitute a
\EM{perfect noncrossing oriented matching} (\nop), whose edges are indicated by grey arcs.}
\label{fig:simplifiedPicture}
\input graphics/simplifiedPicture
\end{figure}

Clearly, we may ``forget'' the \EM{actual} colours (red or green) of the coloured points,
if we remember instead the pattern of radial orientations: The situation is completely
determined by
\bit
\item the geometric positions of the coloured and uncoloured lower and upper points,
	we call this piece of information the \pointconf,
\item together with the pattern of radial orientation; we call this piece of information
	the \pointarr.
\eit
\figref{fig:simplifiedPicture} illustrates this.

Note that for all \pointconf s, the number of coloured \EM{upper} points plus
twice the number of uncoloured \EM{upper} points equals the number of coloured
\EM{lower} points plus twice the number of uncoloured \EM{lower} points.

We call a \pointarr\ \EM{admissible} if
it has the same number of inwardly/outwardly
oriented points.
\EM{Every} admissible \pointarr\ determines together with the
corresponding \pointconf\ a certain colouring (green and red) of lower/upper points,
which represents a certain product of two skew Schur functions
$${\schurf_{\lambda/\mu}}\cdot{\schurf_{\sigma/\tau}},$$
where
\bit
\item $\lambda/\mu$ is determined by the uncoloured points and the green points,
\item $\sigma/\tau$ is determined by the the uncoloured points and the red points.
\eit
Note that there might be no overlay of families of \nilp s that connect
these points (if, for instance, the $i$--th green upper point lies to the left of the
$i$--th green lower point; this would correspond to an $i$--th row of length $<0$
in the corresponding Ferrers diagram: In this case, the corresponding skew Schur function
${\schurf_{\lambda/\mu}}$ is
zero).

But if there \EM{is} an overlay of families of \nilp s that connect
the green and red points, then we may
imagine that the \pointarr\ corresponding to this overlay 
is arranged on a perfect circle:
Observations~\ref{obs:arbitrary-point}, \ref{obs:different-parity} and
\ref{obs:different-orientation} imply that if we draw a straight line connecting two points in the \pointarr\ whenever the corresponding coloured points
are connected by a bicoloured trail, then we will obtain a
\pointmatch\
$\cnr{m}$, i.e.:
\bit
\item \EM{Every} coloured point is incident with an edge in $\cnr{m}$,
\item Every edge in $\cnr{m}$ connects points of \EM{different radial orientation},
\item No two edges in $\cnr{m}$ do \EM{cross} (in the geometric realization as straight lines
	connecting points on a circle).
\eit

In particular, this implies:
\begin{obs}[{\bf Bicoloured trails span equal numbers of inward and outward directed points}]
\label{obs:equal-sizes}
Assume $p$ and $q$ are coloured points in an overlay of \nilp s which are connected by some bicoloured trail. Then the corresponding
\pointarr\ is divided in two parts by the points corresponding to $p$ and $q$:
In each of these
parts, the number of points directed inwards must equal the number of points directed outwards.

Stated otherwise:
Let $\cnr{m}$ be a \pointmatch\ for a \pointarr\ $\cnr{r}$. An arbitrary edge $e$ in $\cnr{m}$ divides 
$\cnr{r}$ into two parts, each of which must contain equal numbers of inward and outward oriented
points.
\end{obs}

We conclude our preparatory considerations with the following observation:
\begin{obs}[{\bf Admissible matchings always can be realized by overlays of \nilp s}]
\label{obs:realizematching}
For every pair $\pas{\cnr{c},\cnr{r}}$, where $\cnr{c}$ is a \pointconf\ and $\cnr{r}$
is an admissible \pointarr\ for $\cnr{c}$, there exists a pair of shapes
$\pas{\lambda/\mu,\sigma/\tau}$ characterized by $\pas{\cnr{c}^\prime,\cnr{r}}$, which has
the \EM{same} upper and lower uncoloured, green and red points as $\pas{\cnr{c},\cnr{r}}$
in the \EM{same} order (i.e., only the \EM{geometric position} of the points may differ),
such that \EM{every} \pointmatch\ $\cnr{m}$ in $\cnr{r}$ can be \EM{realized} by an overlay of \nilp s corresponding to $\pas{\lambda/\mu,\sigma/\tau}$: Maybe the best way to
conceive this is by looking at pictures, 
see \figref{fig:RealizeMatching} and \figref{fig:RealizeMatching2}.
\end{obs}

\begin{figure}
\caption{Illustration: Realizations of matchings. Assume we want to realize a \EM{partial}
\pointmatch\ for some subset $S$ of \EM{coloured} upper points which are \EM{consecutive}
in the corresponding \pointarr. In the corresponding
\pointconf, there might be uncoloured points which are interspersed between the coloured
points. As an example, consider the
partial \pointmatch\ indicated in the upper picture.
Imagine that the red and green paths are
threads dangling down from the respective points: Note that \EM{two} threads are dangling
down from \EM{uncoloured points} (shown as black circles), while only \EM{one} thread is
dangling down from \EM{coloured points} (shown as white circles). Repeat the following step until there
are no more coloured points: 
If two coloured points $p_1$, $p_2$
are connected by an edge $e$ of the partial \pointmatch,
and all points between $p_1$ and $p_2$ (if any) are \EM{uncoloured},
then arrange the threads dangling down from the points $p_1,\dots,p_2$ such that they
alternate in colour and tie together with knots the first and second, the third and fourth, etc.,
of these threads. View these knots as the new uncoloured points and simply \EM{forget}
the points $p_1,\dots,p_2$: Forgotten points are indicated by grey circles. The middle
picture shows an intermediate state, where three edges of the partial \pointmatch\ are
realized by bicoloured trails (shown as thick grey lines). The lower picture shows the final state, where
there are no more coloured points, and all edges of the partial \pointmatch\ are \EM{realized}
by bicoloured trails.
It is clear that the
same algorithm realizes partial \pointmatch s
for subsets of \EM{consecutive coloured} lower points.}
\label{fig:RealizeMatching}
\input graphics/RealizeMatching
\end{figure}

\begin{figure}
\caption{Illustration: Realizations of matchings, 2nd part. Assume that all maximal partial
\pointmatch s on subsets of \EM{consecutive coloured} lower or \EM{consecutive coloured} upper
points already have been
realized, so now we have to realize edges connecting a lower point to an upper point
where there is no coloured point between them (in the right half of the circular arrangement
of points).
The left picture indicates this situation: The edges already realized are shown as black
lines. As in \figref{fig:RealizeMatching}, coloured points are shown as white circles,
uncoloured points as black circles and forgotten points as grey circles.
Obviously, by ``tieing together'' threads of different colours
(as in \figref{fig:RealizeMatching}) combined with ``weaving together'' threads of the same
colour, we can realize the edge (shown as grey arc) that connects these
two coloured points. It is clear that
this algorithm indeed gives ``topological'' bicoloured trails corresponding to the
\pointmatch, and it is also clear, that the ``topological'' situation can be implemented
with lattice paths for appropriately chosen shapes (with large enough distances between
the upper points and lower points).}
\label{fig:RealizeMatching2}
\input graphics/RealizeMatching2
\end{figure}

\section{Applications: Classical and new determinantal identities}\label{sec:apps}
The simple idea for all the identities considered in the rest of this paper
could be stated in an \EM{abstract} manner as follows:
Assume some fixed \pointconf\ $\cnr{c}$.
%
\bit
\item
The set $S$ of \EM{all} overlays of \nilp s corresponding to $\cnr{c}$ and
some fixed admissible \pointarr\ $\cnr{r}$
corresponds to the set of \EM{all} terms in the product
$\schurf_{\lambda/\mu}\cdot\schurf_{\sigma/\tau}$, where the pair of shapes $\pas{\lambda/\mu,\sigma/\tau}$ corresponds to $\pas{\cnr{c},\cnr{r}}$.
\item 
The relation ``corresponds to the same \pointmatch\ as'' is an equivalence relation
on $S$. Stated otherwise, $S$ is \EM{partitioned} into \EM{matching classes} which are determined by the respective \pointmatch s $\cnr{m}_1, \cnr{m}_2,\dots$
of $\cnr{r}$:
$$S=S_{\cnr{m}_1}\cup S_{\cnr{m}_2}\cup \dots, \text{ where }S_{\cnr{m}_i}\cap S_{\cnr{m}_j} = \emptyset\text{ for } i\neq j.$$
(
\figref{fig:InvolStruct} visualizes this concept.)
\item Let $\cnr{m}$ be a fixed \pointmatch\ in $\cnr{r}$, and let $\setof{e_1,\dots,e_k}$ be
a fixed set of edges of $\cnr{m}$: Then the \EM{recolouring} of the bicoloured paths $\setof{b_1,b_2,\dots,b_k}$ corresponding to $\setof{e_1,\dots,e_k}$ in $S_{\cnr{m}}$
effectuates a \EM{weight--preserving} involution
$$
S_{\cnr{m}} \leftrightarrow S^{\prime}_{\cnr{m}}\;,
$$
where  $S^\prime$ is the set of \EM{all} overlays of \nilp s corresponding to $\pas{\cnr{c},\cnr{r}^\prime}$, where $\cnr{r}^\prime$ is the \pointarr\ obtained from $\cnr{r}$ by reversing
the orientation of the points which are incident with the edges 
$\setof{e_1,\dots,e_k}$. Note that except for this reversing of orientations of \EM{points}, the \pointmatch\ $\cnr{m}$ is \EM{unchanged}: Nevertheless, we will call this operation the \EM{reversing of edges} $\setof{e_1,\dots,e_k}$ in $\cnr{m}$.
%
\item If we succeed in ``glueing together'' such \EM{involutions for matching classes} such
that we obtain a \EM{weight preserving bijection}
$$
S_1\cup S_2\cup\dots\leftrightarrow S_1^\prime\cup S_2^\prime\cup \dots,
$$
where
$\setof{S_1,S_2,\dots}$
and $\setof{S_1^\prime,S_2^\prime\dots}$
are two families of sets of overlays
corresponding to the two families
\bit
\item
$\setof{\pas{\lambda_1/\mu_1,\sigma_1/\tau_1}, \pas{\lambda_2/\mu_2,\sigma_2/\tau_2},\dots}$
\item
and $\setof{\pas{\lambda^\prime_1/\mu^\prime_1,\sigma^\prime_1/\tau_1},
\pas{\lambda^\prime_2/\mu^\prime_2,\sigma^\prime_2/\tau^\prime_2},\dots}$
\eit
of pairs of shapes, then this bijection clearly translates to an identity involving sums of products
of (skew) Schur functions
$$
\schurf_{\lambda_1/\mu_1}\cdot\schurf_{\sigma_1/\tau_1} +
\schurf_{\lambda_2/\mu_2}\cdot\schurf_{\sigma_2/\tau_2} + \cdots =
\schurf_{\lambda^\prime_1/\mu^\prime_1}\cdot\schurf_{\sigma^\prime_1/\tau^\prime_1} +
\schurf_{\lambda^\prime_2/\mu^\prime_2}\cdot\schurf_{\sigma^\prime_2/\tau^\prime_2} + \cdots.
$$
\eit
The crucial point is the ``glueing together'' of involutions
for matching classes. Clearly, this amounts to a consistent \EM{rule} which identifies for every
pair $\pas{\cnr{r},\cnr{m}}$ (where $\cnr{r}$ is a \pointarr\ and $\cnr{m}$ is a \pointmatch\ 
in $\cnr{r}$) a \EM{family} $\setof{E_1,\dots,E_k}$ of sets of edges 
to be \EM{reversed}, thus relating $\pas{\cnr{r},\cnr{m}}$ to $\pas{\cnr{r}_1^\prime,\cnr{m}},
\dots,\pas{\cnr{r}_k^\prime,\cnr{m}}$, such that the iteration of this ``edge--reversing
procedure'' yields a ``bipartite substructure'', see \figref{fig:InvolStruct}.
We
call such rule a \EM{recolouring scheme}. In the following, we shall illustrate this
abstract concept by several concrete examples.


\begin{figure}
\caption{Illustration: The picture visualizes the situation where a \EM{recolouring scheme}
amounts to a weight preserving bijection between overlays corresponding to \pointarr s
$\setof{\cnr{r}_1,\cnr{r}_2,\cnr{r}_3}$ and $\setof{\cnr{r}_1^\prime,\cnr{r}_2^\prime}$,
respectively. In the picture it is assumed that there are only 4 \pointmatch s for these
\pointarr s, which are indicated by different shades of grey (white, light grey, grey and
black).
The crucial point is that the recolouring scheme assigns every \pointmatch\ 
to a ``bipartite substructure'' with bipartition classes \EM{of the same size} (the edges of these
``bipartite substructures'' are shown as thin black lines). The simplest case of such ``bipartite substructure'' appears for the black \pointmatch\ 
$\cnr{m}_{\text{black}}$: There is a \EM{mapping} taking $\pas{\cnr{r}_3,\cnr{m_{\text{black}}}}$
to $\pas{\cnr{r}_2^\prime,\cnr{m_{\text{black}}}}$
and vice versa. However, this is not the only possibility: Note that, for example, every white \pointmatch\ 
$\cnr{m}_{\text{white}}$ is related to \EM{two} instances of $\cnr{m}_{\text{white}}$:
The corresponding ``bipartite substructure'' has bipartitions classes 
$\setof{\pas{\cnr{r}_1,\cnr{m_{\text{white}}}},\pas{\cnr{r}_3,\cnr{m_{\text{white}}}}}$ and
$\setof{\pas{\cnr{r}_1^\prime,\cnr{m_{\text{white}}}},\pas{\cnr{r}_2^\prime,\cnr{m_{\text{white}}}}}$ (\EM{both of size $2$)}.}
\label{fig:InvolStruct}
\input graphics/InvolStruct
\end{figure}

\subsection{Dodgson's condensation, revisited and generalized}
\label{sec:Dodgson}
Recall Dodgson's
condensation formula \eqref{eq:Dodgson}, or rather its Schur function equivalent \eqref{eq:DodgsonSchur}.
Given our above preparations,  \figref{fig:Dodgson} contains the
\EM{proof} of \eqref{eq:DodgsonSchur}!

\begin{figure}
\caption{Illustration: Graphical proof of Dodgson's condensation formula. The picture
shows the three \pointarr s that may appear if the bicoloured path starting in the
rightmost coloured upper point $s$ (drawn as black circle) is recoloured. Note that there is only
\EM{one} \pointmatch\ (indicated by grey lines) for the \pointarr s in the left half of the picture, while there are \EM{two} \pointmatch s $\cnr{m}_1,\cnr{m}_2$ for the \pointarr\ in the right half of the picture (the edges incident with $s$ in  $\cnr{m}_1$ and $\cnr{m}_2$ are indicated by grey
lines). The involution effectuated by the recolouring
scheme is indicated by black lines connecting the respective \pointarr s, each of
which represents a product of skew Schur functions.}
\label{fig:Dodgson}
\input graphics/LemmaIllustrationSlide
\end{figure}

In order to see this, consider once again a special case: Let
$\lambda=\pas{\lambda_1,\dots,\lambda_6}=\pas{9, 7, 5, 3, 3, 1}$; for convenience, we
denote $\pas{\lambda_2,\dots,\lambda_5}$ by $\sigma$. Dodgson's formula
starts with the product of Schur functions
$$\schurf_{\lambda}\cdot\schurf_{\sigma} =
\schurf_{\lambda}\cdot\schurf_{\sigma+\constpart{-1}{5}/\constpart{-1}{5}},$$
which by the \lgt--interpretation appears as generating function of overlays of
nonintersecting lattice paths. See \figref{fig:Dodgson2}, where the involutions
visualized in \figref{fig:Dodgson} are ``specialized'' to this concrete example.
These weight--preserving involutions immediately imply \eqref{eq:DodgsonSchur},
and it is obvious that the argument is valid not only for our special example, but
in full generality, whence (by the reasoning following Corollary~\ref{cor:Dodgson})
we have proved \eqref{eq:DodgsonSchur} (and thus \eqref{eq:Dodgson}).
\hfill\qedsymbol

\begin{figure}
\caption{Illustration: Dodgson's condensation formula, by example. The pictures show
the positions of lower and upper points of the overlays corresponding to shapes
derived from $\lambda/\mu=\pas{9, 7, 5, 3, 3, 1}/\allconstpart{0}$ and $\sigma/\tau=\pas{6, 4, 2, 2}/\constpart{-1}{5}$. In the upper
picture,
all
coloured points are green, and the bicoloured path $b$ starting in the rightmost upper
point $s=\pas{8,\n}$ \EM{must} have the rightmost lower point $\pas{-1,1}$ as its
other ending point.  Recolouring $b$ leads to the middle picture. Here, the bicoloured
path starting in $s$ has \EM{two} possibilities: Its other ending point might be $\pas{-1,1}$
again,
but also the leftmost upper point $\pas{-5,\n}$ is possible. For the latter case, recolouring the
bicoloured path connecting $s$ and $\pas{-5,\n}$ leads to the lower picture. Now, as in the
upper picture, the bicoloured path starting in $s$ \EM{must} have $\pas{-5,\n}$ again as
its other ending point. By the above reasoning, this shows immediately the following specialization of
\eqref{eq:DodgsonSchur}:
$$
{\schurf_{\lambda/\mu}}\cdot{\schurf_{\sigma/\tau}}
+
{\schurf_{\lambda^{\prime\prime}/\mu^{\prime\prime}}}\cdot{\schurf_{\sigma^{\prime\prime}/\tau^{\prime\prime}}}
=
{\schurf_{\lambda^{\prime}/\mu^{\prime}}}\cdot{\schurf_{\sigma^{\prime}/\tau^{\prime}}}
$$}
\label{fig:Dodgson2}
\input graphics/Dodgson
\end{figure}

For the proof of Dodgson's condensation formula, we used the following simple
\EM{recolouring scheme}: 
``Fix some nonempty subset $S$ of coloured points
of the same orientation (i.e., all points in $S$ are either all outwards oriented
or all inwards oriented), and in \EM{all} \pointmatch s always reverse all edges incident
with a point in $S$''. We call this rule the \EM{Dodgson recolouring scheme}.

The following Lemma and its proof are a reformulation of \cite[Lemma~15]{fulmek:2001}):
\begin{lem}
\label{lem:DodgsonScheme}
Except for degenerate cases, the Dodgson recolouring scheme always yields a Schur function identity.
\end{lem}
\begin{proof}
Clearly, the Dodson recolouring scheme unambiguously identifies for every \pointmatch\ 
the (single set of) edges to be reversed.
The following proof is merely a clarification of the statement:

Let $\lambda/\mu$ and $\sigma/\tau$ be two skew shapes and consider the \pointconf\ $\cnr{c}$ which corresponds to $\pas{\lambda/\mu,\sigma/\tau}$. Without loss of generality, let $S$
be a nonempty subset of \EM{inwards} oriented points in $\cnr{c}$.

Let $V$ be the set of \EM{all} admissible \pointarr s for $\cnr{c}$ where \EM{all} points
in $S$ have the \EM{same} orientation.
Consider the graph $G$ with vertex $V$, where two vertices
$v_1$, $v_2$ are connected by an edge if and only if there exists a \pointmatch\ $\cnr{m}$
for $v_1$ such that $v_2$ is obtained by applying the Dodgson recolouring scheme (i.e.,
reverse
the edges in $\cnr{m}$ which are incident with a point in $S$).
Clearly, this graph $G$ is \EM{bipartite}: 
$V=I\cup O$, where $I$ is the subset of $V$ with all points in $S$ oriented \EM{inwards}
and $O$ is the subset of $V$ with all points in $S$ oriented \EM{outwards}, and there
is no edge connecting two vertices of $I$ or two vertices of $O$ (see again \figref{fig:Dodgson},
where the bipartite structure is indicated by a dashed line separating the left and right half of the
picture).

Assume that graph $G$ has a connected component $Z$ with at least $2$ vertices
(if there is no such component, we call this a \EM{degenerate case}). Then we
have the following identity for skew Schur functions:
\begin{equation}
\label{eq:lemma-general}
\sum_{\pas{\lambda/\mu,\;\sigma/\tau}\in Z_I}
	\schurf_{\lambda/\mu}\cdot\schurf_{\sigma/\tau}
=
\sum_{\pas{\lambda^\prime/\mu^\prime,\;\sigma^\prime/\tau^\prime}\in Z_O}
	\schurf_{\lambda^\prime/\mu^\prime}\cdot\schurf_{\sigma^\prime/\tau^\prime},
\end{equation}
where $Z_O$ and $Z_I$
denote the sets of pairs of skew shapes corresponding to $\pas{\cnr{c},x}$ for
$x\in O$ and 
$x\in I$, respectively.
\end{proof}

\begin{figure}
\caption{Dodgson's recolouring rule, where the set $S$ consists only of the
rightmost coloured upper point (drawn as black circle), applied to 3 upper/lower points.
}
\label{fig:Dodgson3}
\input graphics/Dodgson3
\end{figure}

\begin{figure}
\caption{Dodgson's recolouring rule, where the set $S$ consists only of the
rightmost coloured upper point (drawn as black circle), applied to 4 upper/lower points.
}
\label{fig:LemmaIllustration}
\input graphics/LemmaIllustration
\end{figure}

We illustrate Lemma~\ref{lem:DodgsonScheme} by the following generalization of
Dodgson's condensation, which might be new. Consider the \pointarr\ with
$k$ upper and $k$ lower points ($k> 1$) which are all green, and always recolour the bicoloured
path starting in the rightmost upper point ($k=2$ corresponds to Dodgson's condensation,
see \figref{fig:Dodgson}; the cases $k=3$ and $k=$ are depicted in \figref{fig:Dodgson2}
and \ref{fig:Dodgson3}).
\begin{thm}
\label{thm:maybenew}
Let $a$ be an $\pas{m+k}\times\pas{m+k}$--matrix, and let $1\leq i_1<i_2<\cdots<i_k\leq m+k$ and
$1\leq j_1<j_2<\cdots<j_k\leq m+k$ be (the indices of) $k$ fixed rows and $k$ fixed columns,
respectively, of $a$. Denote the sets of these (indices of) rows and columns by $R$ and $C$,
respectively.

Let $E\defeq\setof{j_2,j_4,\dots,j_{2\cdot\floor{k/2}}}$ and $O\defeq\setof{i_1,i_3,\dots,i_{2\cdot\ceil{k/2}-1}}$ be the sets of \EM{odd} fixed row indices and of \EM{even} fixed column indices,
respectively. Then we have:
\begin{equation}
\label{eq:maybenew}
\sum_{\substack{S\subseteq E,\\ T\subseteq O,\\ \cardof{S}=\cardof{T}}}
\detof{\deleterowcols{a}{T}{S}}\cdot\detof{\deleterowcols{a}{R\setminus T}{C\setminus S}}
=
\sum_{\substack{S\subseteq E,\\ T\subseteq O,\\ \cardof{S}=\cardof{T}-1}}
\detof{\deleterowcols{a}{T}{S\cup\setof{j_1}}}\cdot\detof{\deleterowcols{a}{R\setminus T}{C\setminus \pas{S\cup\setof{j_1}}}}.
\end{equation}
\end{thm}

\begin{proof}
We shall prove the Schur function identity equivalent to \eqref{eq:maybenew}: Note that
the corresponding \pointconf\ shows $m+k$ lower points, of which $k$ are coloured, and
$m+k$ upper points, of which $k$ are coloured. As always, we concentrate on the coloured
points: Denote the upper coloured points by ${t_1,t_2,\dots,t_k}$,
and the lower coloured points by ${s_1,s_2,\dots,s_k}$
(counted from the right, as always).
Consider the \pointarr\ where \EM{all} coloured points are \EM{green} (see the uppermost
configuration in the left parts of
Figures~\ref{fig:Dodgson2} and \ref{fig:Dodgson3}): This \pointarr\ corresponds to the summand
for $S=T=\emptyset$ in the left--hand side of \eqref{eq:maybenew}.
Clearly, the other end of a bicoloured path starting in the rightmost upper point $t_1$
must either be in the set $O\defeq\setof{s_1,s_3,\dots,s_{2\cdot\ceil{k/2}-1}}$ or in
the set $E\defeq\setof{t_2,t_4,\dots,t_{2\cdot\floor{k/2}}}$. Now observe that it is
possible to recolour points $s_1, t_2, s_3, t_4, \dots$ (in this order;
these recolouring steps are indicated by thick grey lines in Figures~\ref{fig:Dodgson2}
and \ref{fig:Dodgson3}) until
we obtain the \pointarr\ where \EM{all} points in $E\cup O$ are coloured red. (Note that
in this \pointarr, $t_1$ is green if $k$ is even, otherwise $t_1$ is red.)
Now it is easy to see that for an arbitrary choice of subsets $S$, $T$,
where $S\subseteq E\cup\setof{t_1}$, $T\subseteq O$ and $\cardof{S}=\cardof{T}$,
by the recolouring of red upper and lower points in the appropriate order, we can obtain the
situation where the set of red points is precisely the union $S\cup T$: This translates to the assertion.
\end{proof}
As an example, we state the determinantal identity corresponding to
\figref{fig:LemmaIllustration} ($k=4$) for the special case where $a$ is a
$4\times4$--matrix (i.e., $m=0$):
\begin{multline*}
\detof{a} = 
\detprodcminors{a}{1}{1} + \detprodcminors{a}{1,3}{1,2} + \detprodcminors{a}{1,3}{1,4} +
\detprodcminors{a}{3}{1} \\ -
\detprodcminors{a}{1}{2} - \detprodcminors{a}{1,3}{2,4} -
\detprodcminors{a}{3}{2} \\ - \detprodcminors{a}{1}{4} - \detprodcminors{a}{3}{4}.
\end{multline*}

\subsection{(Generalized) Plücker relations}
\label{sec:Pluecker}

There is a particularly simple special case of Lemma~\ref{lem:DodgsonScheme}
(this is a reformulation of
\cite[Lemma~16]{fulmek:2001}):
\begin{cor}
\label{cor:fulmek}
Let $\cnr{r}$ be the \pointarr\ for the
pair of shapes $\pas{\lambda/\mu,\;\sigma/\tau}$, and assume that the orientation of
the points in $\cnr{r}$  is \EM{alternating}.
Consider the set of all \pointarr s which arise by applying Dodgson's recolouring scheme
(for some fixed nonempty set of points of the same orientation) to $\cnr{r}$ and denote the corresponding
set of pairs of
skew shapes by $Q$.

Then we have:
\begin{equation}
\label{eq:lemma-special}
\schurf_{\lambda/\mu}\cdot\schurf_{\sigma/\tau}
=
\sum_{\pas{\lambda^\prime/\mu^\prime,\;\sigma^\prime/\tau^\prime}\in Q}
	\schurf_{\lambda^\prime/\mu^\prime}\cdot\schurf_{\sigma^\prime/\tau^\prime}.
\end{equation}
\end{cor}
\begin{proof}
Consider the \pointarr\ $\cnr{r}^\prime$ corresponding to some Schur function product
$\schurf_{\lambda^\prime/\mu^\prime}\cdot\schurf_{\sigma^\prime/\tau^\prime}$
from the right hand side of \eqref{eq:lemma-special}: By the
combination of Observations \ref{obs:different-orientation} and \ref{obs:different-parity},
it is clear that
re--applying the Dodgson recolouring scheme \EM{must} give the \pointarr\ $\cnr{r}$
corresponding
to the Schur function product $\schurf_{\lambda/\mu}\cdot\schurf_{\sigma/\tau}$.
\end{proof}

\begin{figure}
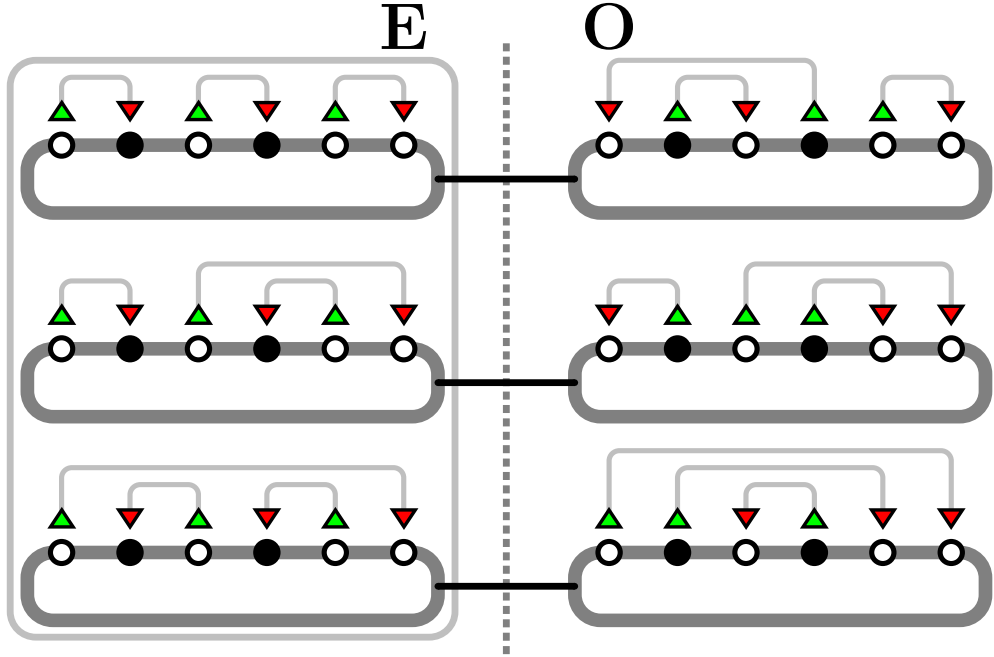

\caption{Illustration of the Pl\"ucker relation \eqref{eq:pluecker} for $m=3$: The left
part shows \EM{three copies} of the \pointarr\ considered in the proof of
Theorem~\ref{thm:Pluecker} with different \pointmatch s; the right part shows the
effect of recolouring the fixed points which are indicated by black circles. Note that
there are \EM{five} possible \pointmatch s in this situation: For the upper and lower
pictures, there are \EM{two} \pointmatch s which yield the same recolouring; these
two \pointmatch s are shown in the left and right parts of the pictures, respectively.}
\label{fig:Pluecker}
\input graphics/Pluecker
\end{figure}

We will use this Corollary for a proof of the Plücker relations (also known as
Grassmann–Pl\"ucker syzygies, see \cite{sturmfels:1993},
or as Sylvester's Theorem, see \cite[section 137]{muir:1933}).
In addition to notation $\range{m}\defeq\setof{1,2,\dots,m}$ we introduce the notation
$$\shiftrange{m}{n}\defeq\setof{n+1,n+2,\dots,n+m}.$$
Moreover, for finite ordered
sets $X\subseteq S$ and $Y$ with $Y\cap S=\emptyset$, $\cardof{Y}=\cardof{X}$, we introduce the
notation $\shuffle{S}{X}{Y}$ for the set $S$ where the elements of $X$ are replaced
by the elements of $Y$ \EM{in the same order}, i.e., if ordered sets $X$ and $Y$ are
given as $X=\pas{x_1,x_2,\dots}$ and $Y=\pas{y_1,y_2,\dots}$, respectively, and $S$ is given as
$$S=\pas{s_1,\dots,s_{(k_1-1)},x_1,s_{(k_1+1)},\dots,s_{(k_2-1)},x_2,s_{(k_2+1)},\dots},$$
then ordered set $\shuffle{S}{X}{Y}$ is given as
$$\shuffle{S}{X}{Y}=\pas{s_1,\dots,s_{(k_1-1)},y_1,s_{(k_1+1)},\dots,s_{(k_2-1)},y_2,s_{(k_2+1)},\dots}.$$
\begin{thm}[Plücker relations]
\label{thm:Pluecker}
Let $a=\pas{a_{i,j}}_{(i,j)\in\range{2m}\times\range{m}}$ be a $2\cdot m\times m$ matrix.
Consider some fixed 
set $R\subseteq\range{m}$.
Then we have
\begin{equation}
\label{eq:pluecker}
\detof{\minor{a}{\range{m}}{\range{m}}}\cdot\detof{\minor{a}{\range{m}}{\shiftrange{m}{m}}}=
\sum_{\substack{S\subseteq\shiftrange{m}{m},\\\cardof{S}=\cardof{R}}}
\detof{\minor{a}{\range{m}}{\shuffle{\range{m}}{R}{S}}}\cdot
\detof{\minor{a}{\range{m}}{\shuffle{\range{2m}\setminus\range{m}}{S}{R}}}.
\end{equation}
\end{thm}
\begin{proof}
We shall prove a Schur function identity which is equivalent to \eqref{eq:pluecker}: Consider
$$
\lambda=\pas{(2\cdot m)\cdot(m-1),(2\cdot m-2)\cdot(m-1),\dots,4\cdot(n-1),2\cdot(n-1)}
$$
and
$$
\sigma=\pas{(2\cdot m-1)\cdot(m-1),(2\cdot m-3)\cdot(m-1),\dots,3\cdot(n-1),(n-1)}
$$
and assume that matrix $\minor{a}{\range{m}}{\range{m}}=\jtmat{\lambda}$ and
matrix $\minor{a}{\shiftrange{m}{m}}{\range{m}}=\jtmat{\sigma}$. It is clear that
the \pointarr\ corresponding to the pair of shapes
$\pas{\lambda/\allconstpart{0},\sigma/\allconstpart{0}}$ has only \EM{upper} coloured points
which alternate in orientation, see the uppermost configuration in the left part of \figref{fig:Pluecker}. Corollary~\ref{cor:fulmek} immediately translates to (the Schur function
equivalent of) \eqref{eq:pluecker}; see \figref{fig:Pluecker} for an illustration.
\end{proof}

Note that the proof (which basically is \EM{contained} in \figref{fig:Pluecker}!)
implies an obvious generalization: Of course, there might be \EM{uncoloured} points
in the \pointconf s associated to the \pointarr s shown in \figref{fig:Pluecker}!
Such uncoloured points amount to a slightly more general Schur function identity
(in \cite[(3.3)]{gurevichPyatovSaponov:2009}, this identity is stated by describing
the situation with certain operations of Ferrers diagrams and proved by using
the Pl\"ucker relations; the connection of this identity to the more general
statements presented here was already explained \EM{adhoc} in \cite{fulmek:2009}),
or to a \EM{common minor} in the products of determinants (as in Theorem~\ref{thm:maybenew}):
\begin{thm}
\label{them:plueckergeneral}
Let $a$ be an $\pas{m+k}\times\pas{m+2\cdot k}$--matrix, and let
$1\leq j_1<j_2<\cdots<j_k\leq m+2\cdot k$ be (the indices of) 
$2\cdot k$ fixed columns of $a$, and set $A\defeq\setof{j_1,\dots,j_k}$ and 
$O\defeq\setof{j_{k+1},\dots,j_{2\cdot k}}$. Let $A^\prime\defeq\range{m+2\cdot k}\setminus A$
and $O^\prime\defeq\range{m+2\cdot k}\setminus O$.
Consider some fixed 
set $R\subseteq A$.
Then we have:
\begin{equation}
\label{eq:plueckergeneral}
\detof{\minor{a}{\range{m+k}}{A^\prime}}
\detof{\minor{a}{\range{m+k}}{O^\prime}}=
\sum_{\substack{S\subseteq O,\\\cardof{S}=\cardof{R}}}
\detof{\minor{a}{\range{m+k}}{\shuffle{A^\prime}{R}{S}}}\cdot
\detof{\minor{a}{\range{m+k}}{\shuffle{O^\prime}{S}{R}}}.
\end{equation}
\end{thm}

\subsection{Laplace expansion}
\label{sec:Laplace}
Assume we are given some \pointarr\ $\cnr{r}$ with at least $m$ lower points $s_1,\dots,s_m$.
We introduce a new recolouring scheme:
Fix lower points $S\defeq\setof{s_{i_1},\dots,s_{i_k}}$ in $\cnr{r}$ ($k\leq m$). For every
\pointmatch\ $\cnr{m}$, reverse  
\bit
\item \EM{one} of the edges of $\cnr{m}$ connecting \EM{two consecutive} coloured upper
	points (we call such edge an upper \EM{handle}), if there is one,
\item else \EM{all} the edges of $\cnr{m}$ starting at a point of $S$.
\eit

We call this the \EM{Laplace recolouring scheme}: Applying this scheme for $S=\setof{j}$ (i.e., $k=1$)
to the
\pointconf\ with $m$ upper and $m$ lower points, $m>1$, all of which are coloured green,
implies the Laplace expansion of the $m\times m$--determinant by its $j$--th column
\begin{equation}
\label{eq:Laplace}
\detof{a} = \sum_{i=1}^m (-1)^{i-j}\cdot\detof{a_{i,j}}\cdot\detof{\deleterowcols{a}{i}{j}}.
\end{equation}
As for Dodgson's condensation, \figref{fig:Laplace} contains the \EM{proof}: The concrete
example presented there illustrates the case $m=4$, $j=1$, and generalizes to a general Schur function identity, which
implies \eqref{eq:Laplace}.\hfill\qedsymbol
\begin{figure}
\caption{The proof of the Laplace expansion of the $4\times 4$--determinant is contained
in the pictures below. The starting point $s$ of the bicoloured trail to be recoloured is
indicated by a coloured circle,
and the possible connections by bicoloured trails are indicated by grey arcs. Note that
for some partition $\lambda=\pas{\lambda_1,\dots,\lambda_4}$, 
the involutions
(indicated in the pictures by thick black lines) imply that the generating function of the
left half,
$$
\schurf_{\addtoparts{\lambda}{1}/{\constpart{1}{4}}}\cdot\schurf_{\pas0} +
\schurf_{\deletepart{\addtoparts{\lambda}{1}}{2}}\cdot\schurf_{\pas{\lambda_2-1}}+
\schurf_{\deletepart{\addtoparts{\lambda}{1}}{4}}\cdot\schurf_{\pas{\lambda_4-3}},
$$
is equal to the generating function of the right half,
$$
\schurf_{\addtoparts{\lambda}{1}/{\constpart{1}{4}}}\cdot\schurf_{\pas0} +
\schurf_{\deletepart{\addtoparts{\lambda}{1}}{2}}\cdot\schurf_{\pas{\lambda_2-1}}+
\schurf_{\deletepart{\addtoparts{\lambda}{1}}{4}}\cdot\schurf_{\pas{\lambda_4-3}}.
$$
By the connections between minors of $\jtmat{\sigma/\tau}$ and operations on partitions
$\sigma$ and $\tau$ (which were explained before Corollary~\ref{cor:Dodgson}), this
means $\detof{\jtmat\lambda}=
\detof{\deleterowcols{\jtmat\lambda}{1}{1}}\cdot\pas{\jtmat\lambda}_{1,1}-
\detof{\deleterowcols{\jtmat\lambda}{2}{1}}\cdot\pas{\jtmat\lambda}_{2,1}+
\detof{\deleterowcols{\jtmat\lambda}{3}{1}}\cdot\pas{\jtmat\lambda}_{3,1}-
\detof{\deleterowcols{\jtmat\lambda}{4}{1}}\cdot\pas{\jtmat\lambda}_{4,1},
$
which implies the Laplace expansion (by the first column).
}
\label{fig:Laplace}
\input graphics/Laplace
\end{figure}

For the generalization of Laplace's expansion (see \cite[section 93]{muir:1933}),
we introduce the following notation: Let $X=\setof{x_1,\dots,x_m}\subset\Z$ be an ordered
set, and let $S=\setof{x_{i_1},\dots,x_{i_k}}$, $k\leq m$, be a subset. In this
situation, we define
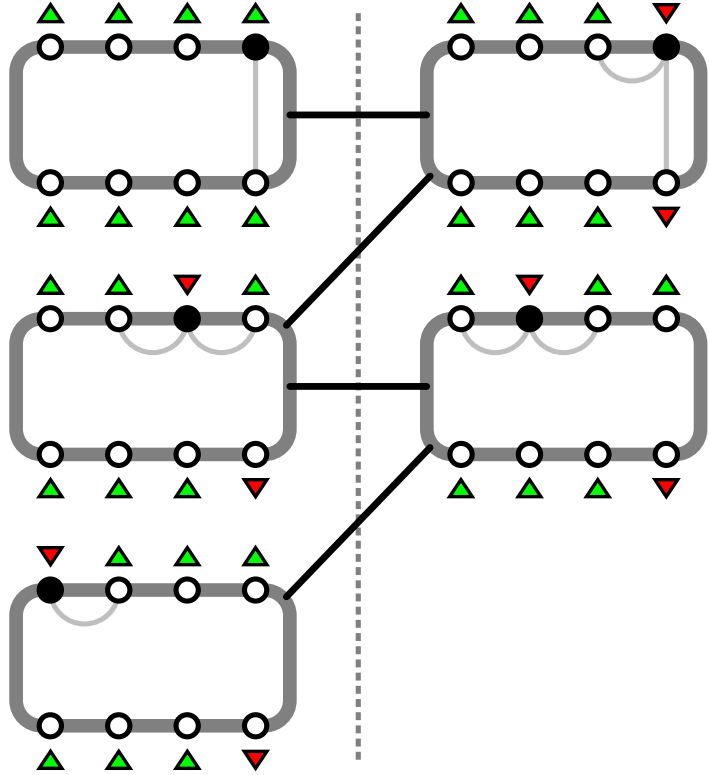
$$\sumset{S}{X}\defeq\sum_{j=1}^k i_j.$$
\begin{thm}[Laplace's Theorem]
\label{thm:Laplace}
Let $a=\pas{a_{i,j}}_{(i,j)\in\range{m}\times\range{m}}$ be an $m\times m$--matrix.
Consider some fixed set $I\subseteq\range{m}$.
Then we have
\begin{equation}
\label{eq:laplace-general}
\detof{a} = 
\sum_{\substack{J\subseteq\range{m},\\\cardof{J}=\cardof{I}}}
\pas{-1}^{\sumset{I}{\range{m}}+\sumset{J}{\range{m}}}\cdot
\detprodcminors{a}{I}{J}.
\end{equation}
\end{thm}

\begin{figure}
\caption{Illustration of Laplace's Theorem, for $m=5$ and $J=\setof{2,4}$: The lower points
corresponding to $J$ are indicated by black circles. The \pointmatch\ shown in the picture
has two upper handles
(edges connecting neighbouring upper points); handles that have been ``reversed'' (i.e., their
endpoints were recoloured) are drawn in black. The important point is that the same number
of pairs (\pointarr,\pointmatch) appears in the two bipartition classes (corresponding to
the left and the right half of the picture); see also \figref{fig:InvolStruct} where
this is situation appears in a more ``abstract'' way.}
\label{fig:Laplace2}
\input graphics/Laplace2
\end{figure}


\begin{proof} As always, we consider the equivalent Schur function identity:
The \pointarr\ corresponding to the left--hand side of \eqref{eq:laplace-general}
consists of
\bit
\item $m$ upper points, which are outward oriented,
\item and $m$ lower points, which are inward oriented
\eit
(stated otherwise: All $2\cdot m$ coloured points are green; see \figref{fig:Laplace2} for an
illustration). We may assume $0<\cardof{I}<m$ (otherwise, there is nothing to prove),
i.e., $I$ corresponds to $k$ lower points $s_{i_1},\dots,s_{i_k}$, $0<k<m$.

Now we apply the Laplace recolouring scheme.
Note that there \EM{always} is an upper handle if there is an edge connecting two
upper points.
If there is more than one upper handle, the recolouring scheme does \EM{not}
describe a \EM{mapping}, but a \EM{multi--valued relation} on pairs $\pas{\cnr{r},\cnr{m}}$, where
$\cnr{r}$ is a \pointarr\ and $\cnr{m}$ is a \pointmatch\ in $\cnr{r}$: It is easy
to see that nevertheless a \EM{bipartite substructure} occurs (see \figref{fig:InvolStruct}): Starting with some arbitrary but fixed pair
$\pas{\cnr{r},\cnr{m}}$, by the Laplace recolouring scheme all the pairs $\pas{\cnr{r}^{\prime\prime},\cnr{m}}$ which
are obtained by reversing an \EM{even} number of upper handles, end up in the same bipartition
class as $\pas{\cnr{r},\cnr{m}}$: Clearly, these are of the same number as all the pairs $\pas{\cnr{r}^{\prime},\cnr{m}}$ which
are obtained by reversing an \EM{odd} number of upper handles (constituting the other bipartition
class).
\end{proof}

Note that the same proof also works for \pointconf s which contain \EM{uncoloured}
points: This amounts to a slightly more general statement (see \cite[section 148]{muir:1933})
\begin{thm}
\label{thm:muir148}
Let $a$ be an $\pas{m+k}\times\pas{m+ k}$--matrix, and let
$1\leq i_1<i_2<\cdots<i_m\leq m+k$ and
$1\leq j_1<j_2<\cdots<j_m\leq m+k$ be (the indices of) $k$ fixed rows and
$k$ fixed columns of $a$. Denote the set of these (indices of) rows and columns by $R$ and $C$,
respectively.
Consider some fixed 
set $I\subseteq R$.
Then we have:
\begin{equation}
\label{eq:muir148}
\detof{a}\cdot\detof{\deleterowcols{a}{R}{C}} = 
\sum_{\substack{J\subseteq C,\\\cardof{J}=\cardof{I}}}
\pas{-1}^{\sum{I}{R}+\sumset{J}{C}}\cdot
\detprodcminors{a}{R}{S}.
\end{equation}
\end{thm}

\bibliography{paper}

\end{document}